\documentclass[11pt,a4paper]{article}

\usepackage[english]{babel}
\usepackage[nottoc]{tocbibind}
\usepackage[utf8]{inputenc}
\usepackage{amsmath}
\usepackage{graphicx}
\usepackage{epstopdf}
\usepackage{amssymb, ulem}
\usepackage{color}
\usepackage[dvipsnames]{xcolor}
\usepackage{float}
\newtheorem{theorem}{Theorem}[section]
\newtheorem{lemma}[theorem]{Lemma}
\newtheorem{assumption}{Assumption}

\newtheorem{remark}[theorem]{Remark}

\newtheorem{example}[theorem]{Example}
\catcode`@=11 \@addtoreset{equation}{section}
\renewcommand\theequation{\thesection.\@arabic\c@equation}
\def\eop{{\ \vrule height 7pt width 7pt depth 0pt}}
\makeatother

\usepackage{nicefrac}

\newcommand{\sfrei}{}

\newcommand{\vertiii}[1]{{\left\vert\kern-0.25ex\left\vert\kern-0.25ex\left\vert #1 
		\right\vert\kern-0.25ex\right\vert\kern-0.25ex\right\vert}}

\usepackage{float}

\title{\bf \Large An implicitly extended Crank-Nicolson scheme for the heat equation on a time-dependent domain}

\author{Stefan Frei~\thanks{Department of Mathematics and Statistics, University of Konstanz, 78457 Konstanz, Germany, stefan.frei@uni-konstanz.de (Corresponding author)},\, Maneesh Kumar Singh~\thanks{Department of Mathematics, Imperial College London, SW7 2AZ, London, UK, maneesh-kumar.singh@imperial.ac.uk}}

\date{}

\begin{document}

\maketitle


\begin{abstract}
We consider a time-stepping scheme of Crank-Nicolson type for the heat equation on a moving domain in Eulerian coordinates. As the spatial domain varies between subsequent time steps, an extension of the solution from the previous time step is required. Following Lehrenfeld \& Olskanskii [ESAIM: M2AN, 53(2):\,585-614, 2019], we apply an implicit extension based on so-called ghost-penalty terms. For spatial discretisation, a cut finite element method is used. We derive a complete a priori error analysis in space and time, which shows in particular second-order convergence in time under a parabolic CFL condition. Finally, we present numerical results in two and three space dimensions that confirm the analytical estimates, even for much larger time steps.
\end{abstract}

\section{Introduction}

Partial differential equations (PDEs) posed on moving domains are significant in many areas of science and engineering. They arise for example in flow problems around moving structures, such as pumps~\cite{Becker1995}, wind or water turbines~\cite{Porteetal2020}, within moving objects~\cite{caucha2018finite}, or as sub-problems in fluid-structure interactions or multiphase flows. Fluid-structure interactions arise in aerodynamical applications like flow around airplanes or parachutes~\cite{Steinetal2000}, in biomedical problems such as blood flow through the cardiovascular system~\cite{Peskin1972, VandeVosseetal2003, Formaggiaetal2010} or the airflow within the respiratory system~\cite{WallRabczuk2008} and even in tribological applications~\cite{KnaufFreiRichterRannacher}. Multiphase problems include for instance
gas-liquid and particle-laden gas flows~\cite{GurrisKuzminTurek, crowe1997, JudakovaBause}, rising bubbles~\cite{Hysingetal2009}, droplets in microfluidic devices~\cite{ClausKerfriden2019} or the simulation of tumor growth~\cite{Garckeetal2018}. For further details and applications we refer to the textbooks~\cite{richter2017fluid, bazilevs2013computational} and~\cite{GrossReusken}, respectively.

In this article we consider the time discretisation of a parabolic model problem (namely the heat equation) which is posed on a moving domain $\Omega(t) \subset \mathbb{R}^{d}\,(d = 2, 3)$ that evolves smoothly in time for $t \in I = [0, t_{max}]$:
\begin{equation}\label{C1}
\begin{array}{ll}
  u_{t}-\Delta u =f \quad \mbox{in} \,\,\Omega(t),
   \qquad  u=0 \quad \mbox{on} \,\,\partial \Omega(t),  
    \qquad u(x,0)=u_{0}(x) \quad \mbox{in} \,\,\Omega(0).
\end{array}
\end{equation}


In literature, two major numerical approaches can be found for the simulation of partial differential equations on moving domains: the Arbitrary Lagrangian Eulerian (ALE) \cite{donea1982arbitrary,donea2004arbitrary} approach, where the equations are transformed to an arbitrary reference domain  which is independent of time, and Eulerian approaches, where the equations are solved in the time-dependent Eulerian framework~\cite{dunne2006eulerian, FreiDiss, Lc_esiam19}. 

The ALE approach is a popular technique for the numerical simulation of PDEs on moving domains, in particular for flow problems~\cite{hughes1981lagrangian,hirt1997arbitrary}. For details, we refer to the 
textbooks~\cite{bazilevs2013computational,richter2017fluid} and reference cited therein.
Convection-diffusion problems on moving domains were, for example, solved in~\cite{ganesan2017ale,srivastava2020local} using a stabilised ALE method. The ALE approach is very attractive in the case of moderate domain movements, but shows problems when the shape of the domain changes significantly in time. In particular, topology changes of $\Omega(t)$, as occurring for example in contact problems or considering the separation or union of bubbles can not be modelled by means of an ALE approach~\cite{richter2017fluid,ClausKerfriden2019,burman2020nitsche,burman2021mechanically}. Another example of extreme variations of $\Omega(t)$ are so-called fingering phenomena, which can be frequently observed in multi-phase flows or even for tumor growth~\cite{Garckeetal2018}. 

In such cases, a numerical approach that discretises the equations directly in the moving Eulerian coordinate framework is preferable. The Eulerian framework is also the coordinate framework, which is typically used to model flow problems and consequently, in multi-phase flows~\cite{GrossReusken} and fluid-structure interactions with large displacements~\cite{dunne2006eulerian, BurmanFernandez2014, FreiDiss}. However, as the domains $\Omega(t)$ to be discretised vary with time $t$, additional difficulties arise concerning a proper and accurate discretisation, both in space and in time.  


In recent decades, a great amount of works have been contributed concerning the spatial discretisation of curved or moving boundaries by means of finite elements. The techniques can be categorised in \textit{fitted} and \textit{unfitted} finite element methods. In \textit{fitted} methods, 
the boundary $\partial\Omega(t)$ is resolved in each time step by the finite element mesh~\cite{FeistauerSobotikova1990, BrambleKing1996, frei2014locally}. If the domain is time-dependent, this means that new meshes need to be created in each time step. Several approaches have been proposed to alleviate this issue, such as the locally fitted finite element method~\cite{frei2014locally,frei2020locally}, which is based on a fixed coarse and a variable fine mesh. However, different issues might arise, such as anisotropic fine cells that complicate the numerical discretisation~\cite{FreiPressure}. 

The idea in \textit{unfitted} finite element methods, on the other hand, is to use the same finite element mesh for all times $t$, independently of the position of the boundary $\partial\Omega(t)$. A popular approach is the cut finite element method (CutFEM) \cite{burman_ijnme15,burman2016cut,hansbo2014cut,massing2018stabilized,hansbo2002unfitted,zahedi2017space}, 
where cells of the finite element mesh are cut into parts that lie inside $\Omega(t)$ and parts outside for numerical integration. Boundary values are then incorporated weakly by means of Nitsche's method~\cite{Nitsche70}. The method shows similarities to the extended finite element method \cite{daux2000arbitrary,chessa2003extended,fries2010extended} and the generalised finite element method~\cite{BabuskaBanarjeeOsborn2004}, where the finite element spaces are enriched by suitable functions to account for the position of the boundary.

Much less works can be found in literature concerning a proper time discretisation on moving domains. In the case of moving domains, standard time discretisation based on the method of lines is not directly applicable. The reason is that the domain of definition of the variables changes from time step to time step.
As an example consider the finite difference discretisation of the time derivative within a variational formulation ($\Delta t = t_n - t_{n-1}$)
\begin{align*}
 (\partial_t u_h(t_n), \phi_h^n)_{\Omega(t_n)} \approx \frac{1}{\Delta t} (u_h(t_n)- u_h(t_{n-1}), \phi_h^n)_{\Omega(t_n)}.
\end{align*}
The function $u_h(t_{n-1})$ is only well-defined on $\Omega(t_{n-1})$, but is needed on $\Omega(t_n)$.

A possible remedy is to use characteristic-based approaches based on trajectories that follow the motion of the domain, see e.g.~\cite{HechtPironneau2016}. Similar time-stepping schemes result when applying the ALE method only locally within 
one time step and projecting back to the original reference frame after each step~\cite{Codina2009}, or based on Galerkin time discretisations with modified Galerkin spaces~\cite{FreiRichter17}.
The disadvantage of these approaches is the necessity for a  
projection that needs to be computed within each or after a certain number of steps.

A further alternative are space-time approaches~\cite{HansboLarsonZahedi2016, Lehrenfeld2015}, where a $d+1$-dimensional domain is discretised. These are, however, computationally demanding, in particular within complex three-dimensional applications. 
The implementation of higher-dimensional discretisations and accurate quadrature formulas pose additional challenges. If a discontinuous Galerkin approach is applied in time for the test functions, the formulation decouples in certain time intervals and can be seen as an Eulerian time-stepping scheme~\cite{HansboLarsonZahedi2016, zahedi2017space, FrachonZahedi2019}.

In this work, we follow a slightly different approach first used by Schott~\cite{schott2017stabilized} and later analysed by
Lehrenfeld \& Olshanskii~\cite{Lc_esiam19}. Here, the idea is to define
extensions of the solution $u(t_{n-1})$ from previous time steps to a domain $\Omega_\delta(t_{n-1})$ that spans at least $\Omega(t_n)$. On the finite element level these extensions can be incorporated implicitly in the 
time-stepping scheme by so-called \textit{ghost penalty} stabilisations~\cite{Burman2010} to a sufficiently large domain. These techniques have originally been proposed to extend 
the coercivity of elliptic bilinear forms from the physical to the computational domain in the context of CutFEM or fictitious domain approaches~\cite{Burman2010}.

Lehrenfeld \& Olshanskii \cite{Lc_esiam19} analysed the so extended Backward Euler method in detail for a convection-diffusion problem and gave hints on how to transfer the argumentation to the second-order backward difference scheme (BDF2). Recently, the analysis has been extended to higher order in space and time using an isoparametric finite element approach~\cite{lou2021isoparametric}. In \cite{burman2022eulerian,vonunfitted21}, extended BDF time-stepping schemes were applied and analysed for the non-stationary Stokes equations on moving domains. 

The reason why only BDF-type time-stepping schemes have been considered in previous works, is that in these schemes spatial derivatives appear only on the {\sfrei ''new'' time step, i.e.\,$\nabla u(t_{n})$}. We will see below that the appearance of additional derivatives on $u(t_{n-1})$ will complicate the error analysis severely. This paper gives a first step towards the analysis of time-stepping schemes that require derivatives at different time instants, such as the Crank-Nicolson method, the Fractional-step-$\theta$-, implicit Runge-Kutta- or Adams-Bashforth schemes.

As a first step, we focus in this work on the popular Crank-Nicolson time-stepping scheme. 
Up to now, it has been largely open, if and under what conditions a Crank-Nicolson-type scheme can be used within an Eulerian time discretisation on moving domains.
We give a detailed stability and convergence analysis. {\sfrei While the analysis requires a strong parabolic CFL condition of type $\Delta t\leq ch^2$, our numerical results indicate that the scheme is stable also for much larger time steps.}


The article is organised as follows: In Section \ref{Sec:2}, we introduce the discretisation of the model problem (\ref{C1}) in time and space. Section \ref{Sec:3} presents a stability analysis for the fully discrete scheme using a CFL condition. In Section \ref{Sec:4}, we show a detailed a priori convergence analysis. Numerical experiments in two and three space dimensions are presented in Section \ref{Sec:5}. Section \ref{Sec:6} summarises this  article with some concluding remarks.

\section{Discretisation}\label{Sec:2}
In this section, we present the numerical approximation of the model problem (\ref{C1}). We start with discretisation in time, and continue with the spatial discretisation of the resulting time-discrete formulation.

\subsection{Temporal discretisation}
For time discretisation, we divide the time interval of interest $I=[0,t_{max}]$ in intervals $I_{n}=(t_{n-1},t_{n}]$. For simplicity, we take a uniform time step $\Delta t = \dfrac{t_{max}}{N}$ and define  $t_{n}=n \Delta t$. We define the domain $\Omega^{n}:=\Omega(t_{n})$ with boundary $\Gamma^{n}:=\Gamma(t_{n})$ and write ${u^{n}}=u(t_n)$ for the exact solution of the continuous problem (\ref{C1}) at time $t_n$.

\begin{figure}[t]
\centerline{
\resizebox*{7cm}{!}{\includegraphics{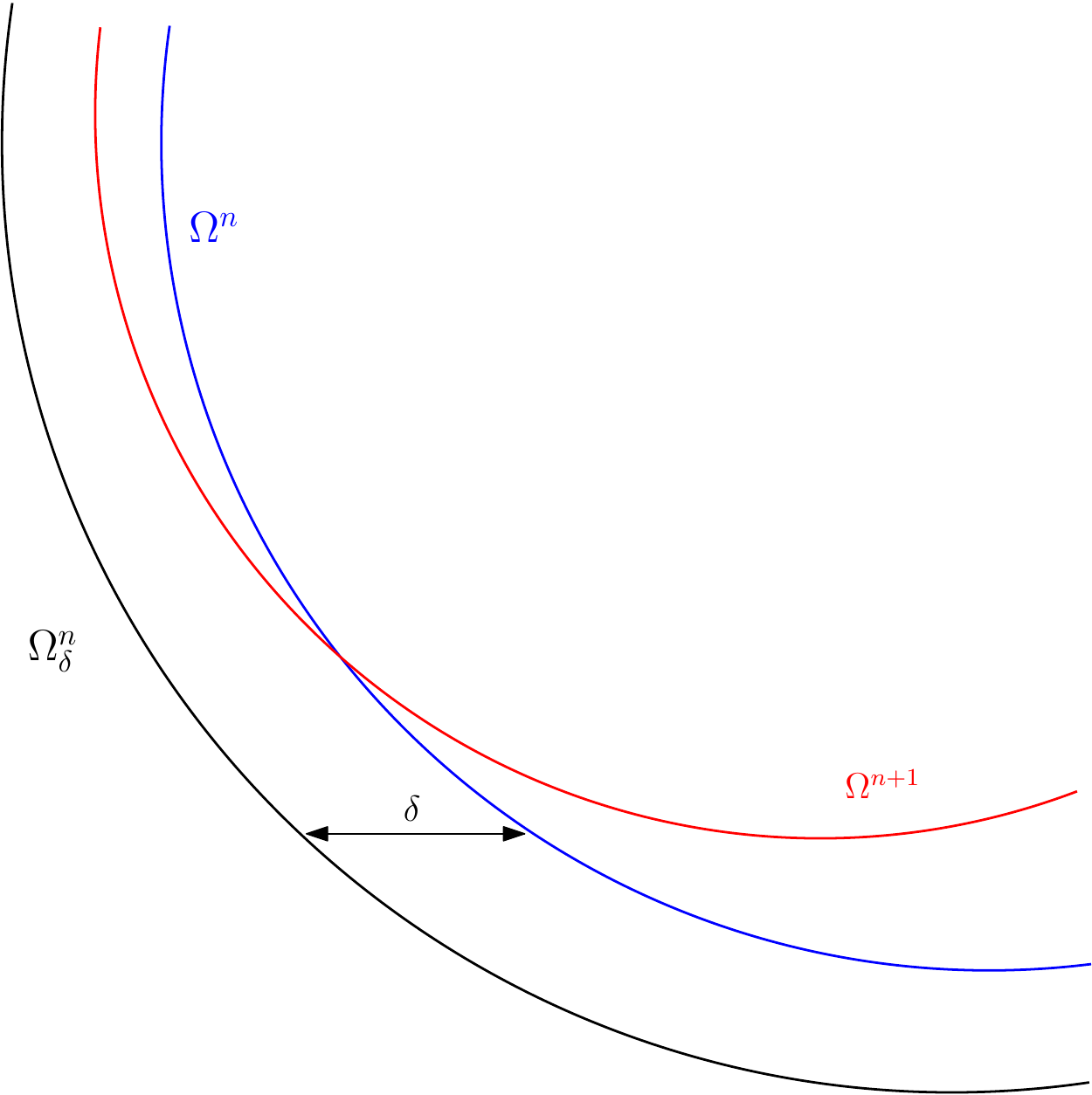}}%
} 
\caption{{\label{fig01} Illustration of the domains $\Omega^{n}$, $\Omega^{n+1}$ and the extension $\Omega^{n}_{\delta}$.}}
\end{figure}

A $\delta$-neighborhood of $\Omega(t)$ at time step $n$ is chosen large enough such that 
$(\Omega^n \cup \Omega^{n+1}) \subset \Omega^{n}_{\delta}$, see Figure~\ref{fig01}. Therefore we choose
\[
\delta \geq w_{max}\Delta t, \qquad w_{max} = \sup_{t \in I, \, x \in \partial \Omega(0)}\|\partial_{t}T(x,t) \cdot n\|.
\]
The required regularity of the domain mapping $T$ will be ensured in Assumption~\ref{ass:regularity} below. For the error analysis we will also assume the upper bound
\begin{align}\label{DeltaUpperBound}
\delta \leq c_{\delta} w_{max}\Delta t
\end{align}
with a constant $c_\delta >1$. Finally, we introduce the following notations for some space-time domains
\begin{align*}
Q:=\mathop{\cup}_{t \in I}\{t\} \times \Omega(t), \quad Q^{n}:=\mathop{\cup}_{t \in I_n}\{t\} \times \Omega(t)
\quad Q^{n}_{\delta}:=\mathop{\cup}_{t \in I_n}\{t\} \times \Omega_{\delta}(t),\quad \hat{Q} = \Omega(0) \times [0,t_{max}].
\end{align*}
Now, the Crank-Nicolson method applied to~\eqref{C1} writes formally
\begin{equation}\label{tc1}
    \dfrac{{u^{n}}-{u^{n-1}}}{\Delta t}-\dfrac{1}{2}(\Delta {u^{n}}+\Delta {u^{n-1}})=\dfrac{1}{2}({f^{n}}+{f^{n-1}}),\quad x \in \Omega^{n}. 
\end{equation}
The main issue of this formulation is that ${u^{n-1}}$ is needed on $\Omega^n$, while it is defined on $\Omega^{n-1}$. Thus, we will add implicit extension operators below to define ${u^{n}}$ on $\Omega_\delta^n \supset \Omega^{n+1}$, where it is needed in the following time step. Similarly, ${f^{n-1}}$ might be undefined on $\Omega^n\setminus \Omega^{n-1}$. If ${f^{n-1}}$ is given analytically, it can typically be extended in a canonical way to $\Omega^n$. To cover different scenarios, we do not want to restrict the analysis in this work to a particular extension, but assume only that ${f^{n-1}}$ is smoothly extended to $\Omega^n$.

In this article, we will use the abbreviation $c$ to refer to a generic positive constant, which is independent of discretisation parameters ($\Delta t, \,h$) and the relative positions of the boundary with respect to the mesh.

\subsubsection{Extension operator}

In this part, we introduce an extension operator of the exact solution $u$ to extend variables to larger domains, as the spatial domain evolves. We make the following assumption (see also \cite[Assumption 3.2]{burman2022eulerian}) for the analysis of this article.

\begin{assumption}\label{ass:regularity}
The boundary of the initial domain $\Omega(0)$ is assumed to be piecewise smooth
and Lipschitz, and the domain motion $T(t)$ is a $W^{1,\infty}$-diffeomorphism for each t, that fulfills $T\in W^{r,\infty}(\hat{Q})$, where $r=\max\{3,m+1\}$ and $m$ is the polynomial degree of the finite element space defined in the following subsection.
\end{assumption}

By using assumption~\ref{ass:regularity}, there exist $W^{r,\infty}$-stable extension operators $E^{n}$ from $\Omega^n$ to $\Omega_{\delta}^n$ that satisfy the following analytical properties:
\begin{align}
    \|E^{n}u-u\|_{W^{m+1,p}(\Omega)}=0, \quad \|E^{n}u\|_{W^{m+1,p}(\Omega^{n}_{\delta})} \leq c \|u\|_{W^{m+1,p}(\Omega^{n})}, \label{exteq1}\\[6pt]
    \|\partial_{t}E^{n}u\|_{H^{m}(\Omega^{n}_{\delta})} \leq c \left( \|u\|_{H^{m+1}(\Omega^{n})}+\|\partial_{t} u\|_{H^{m}(\Omega^{n})}\right), \label{exteq2}\\[6pt]
    \|\partial^{3}_{t}E^{n}u\|_{L^{\infty}(Q^{n}_{\delta})} \leq c \|u\|_{W^{3,\infty}(Q)}. \label{exteq3}
\end{align}
The properties (\ref{exteq1}) and (\ref{exteq2}) are discussed in~\cite{burman2022eulerian}. In an analogous way, one can derive the estimate for the third-order time derivative in (\ref{exteq3}).

\subsection{Spatial discretisation}

For spatial discretisation, we introduce a polygonal domain $D$, which is chosen large enough, such that $\Omega_{\delta}^{n} \subset D$ for all $n$. We introduce a quasi-uniform family of triangulations $({\cal T}_h)_{h>0}$ of $D$ with maximum cell size $h$, which will serve as background meshes.  

In each time step, we extract from ${\cal T}_h$ all cells of non-empty intersection with $\Omega_\delta^n$ and define
\begin{align*}
\mathcal{T}^{n}_{h,\delta}:= \{ K\in {\cal T}_h: \, K\cap \Omega_\delta^n \neq \emptyset\}.
\end{align*}
We write $\Omega^{n}_{h,\delta}$ for the domain spanned by all cells $K \in \mathcal{T}^{n}_{h,\delta}$ and 
define the following finite element space: 
\[
V^{n,m}_{h}:=\{v \in C(\Omega^{n}_{h,\delta}), v|_{K} \in P_{m}(K) \,\,\forall K \in \mathcal{T}^{n}_{h,\delta}\},\quad m \geq 1.
\]

The set of elements that lie (at least partially) outside of $\Omega^{n-1}$, but in $\Omega^n$, will be of particular interest in the analysis. The domain spanned by them will be denoted by 
\begin{align*}
\mathcal{S}^{n,n-1}_{h} :=\mathop{\cup}_{K \in \mathcal{T}^{n, n-1}_{h}} K, \qquad \text{where } \mathcal{T}^{n, n-1}_{h} := \left\{ K\in {\cal T}_{h,\delta}^n, \, K \cap (\Omega^{n} \setminus \Omega^{n-1}) \neq \emptyset\right\}.
\end{align*}

\begin{figure}[t]
\centerline{
\begin{tabular}{cc}
\resizebox*{9cm}{!}{\includegraphics{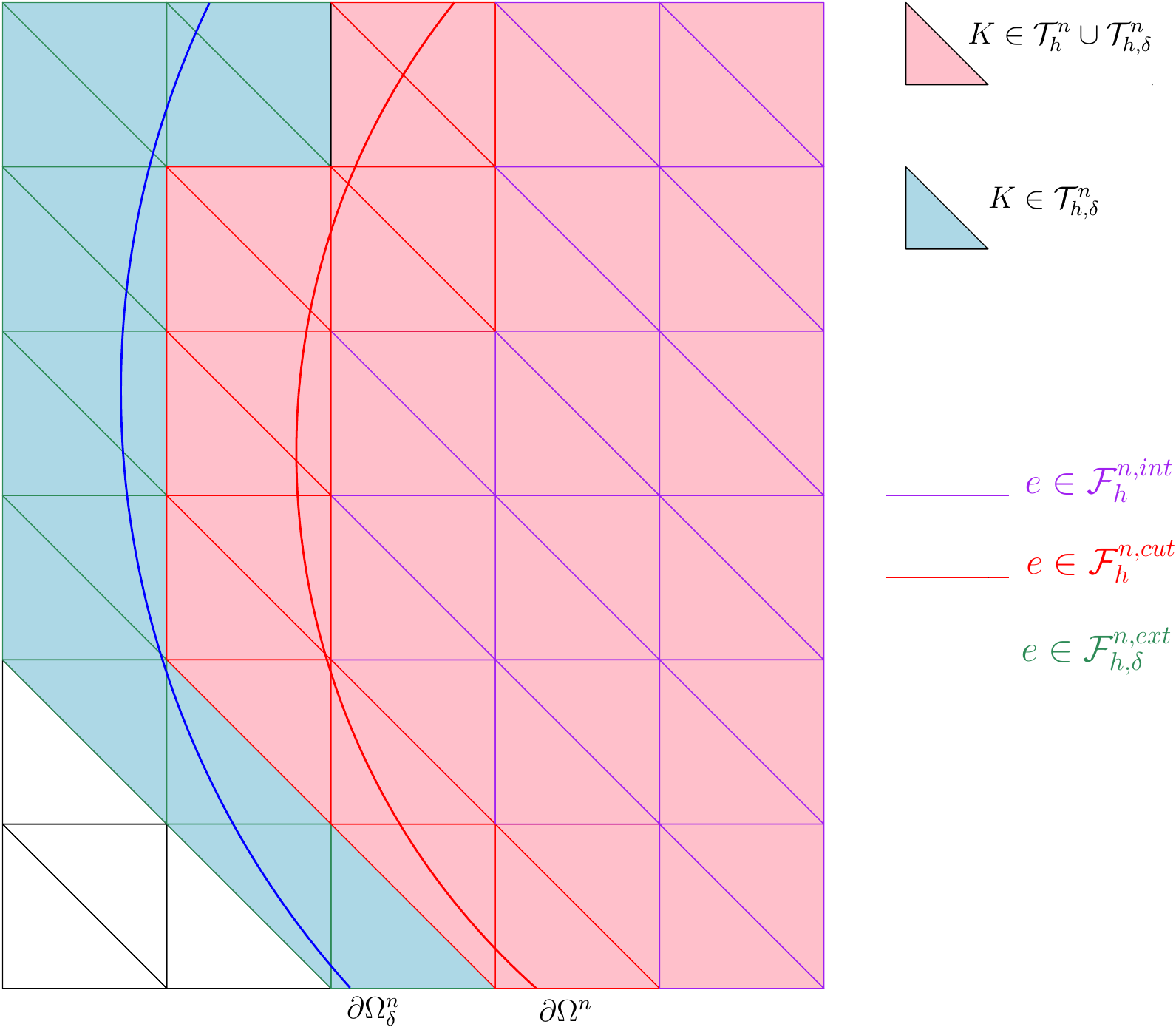}}%
& 
\resizebox*{7cm}{!}{\includegraphics{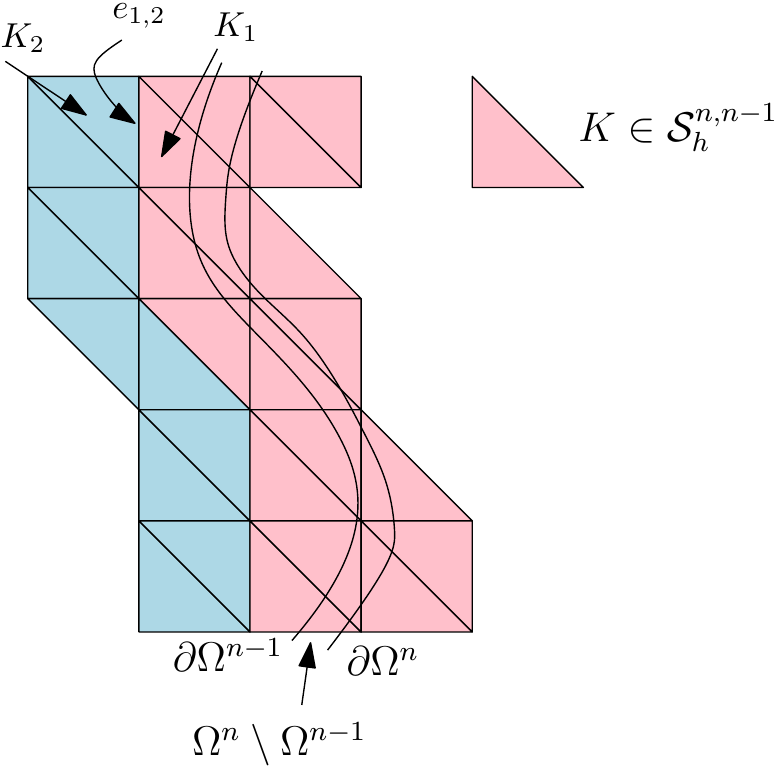}}%
\end{tabular}
} 
\caption{\label{fig:domaindisp} 
\textbf{Left}: Illustration of the triangulations $\mathcal{T}^{n}_{h}$ and $\mathcal{T}^{n}_{h,\delta}$ and the sets of facets $\mathcal{F}^{n}_{h,\delta} = \mathcal{F}^{n,int}_{h,\delta}\cup \mathcal{F}^{n,cut}_{h,\delta} \cup \mathcal{F}^{n,ext}_{h,\delta}$.  \textbf{Right}: Set of boundary cells $\mathcal{S}^{n,n-1}_{h}$.}
\end{figure}

Moreover, we introduce the following notations for the facets of $\mathcal{T}^{n}_{h,\delta}$, see Figure~\ref{fig:domaindisp} for an illustration:
\begin{itemize}
\item 
$\mathcal{F}^{n}_{h,\delta}$ : the set of interior facets of $\mathcal{T}^{n}_{h,\delta}$.

\item 
$\mathcal{F}^{n,int}_{h,\delta}$ : the set of facets that belong exclusively to elements $K \in \mathcal{T}^{n}_{h,\delta}$
that lie completely in the interior of $\Omega^{n}$.

\item
$\mathcal{F}^{n,cut}_{h,\delta}$ : the set of facets that belong to some element $K \in \mathcal{T}^{n}_{h,\delta}$ with $K \cap \partial \Omega^{n} \neq \emptyset$.

\item 
$\mathcal{F}^{n,ext}_{h,\delta}$ :  the set of remaining facets of $\mathcal{F}^{n}_{h,\delta}$, i.e.~$\mathcal{F}^{n,ext}_{h,\delta}:= \mathcal{F}^{n}_{h,\delta} \setminus \left( \mathcal{F}^{n,int}_{h,\delta} \cup \mathcal{F}^{n,cut}_{h,\delta} \right)$

\item 
$\mathcal{F}^{n,g}_{h,\delta}$ := $\mathcal{F}^{n,cut}_{h,\delta} \cup \mathcal{F}^{n,ext}_{h,\delta}$.
\end{itemize}


\begin{assumption}[CFL condition]\label{ass:CFL}
We assume the parabolic CFL condition $\Delta t\leq c_{CFL} h^2$, where $c_{CFL}$ is an arbitrary constant for $m=1$, while we assume $c_{CFL}$ sufficiently small for $m>1$.
\end{assumption}


\noindent The inequality~\eqref{DeltaUpperBound} and the CFL condition (Assumption~\ref{ass:CFL}) lead to 
\begin{align}\label{deltah2}
\delta \leq c_{\delta} w_{max}\Delta t \leq c  w_{max} h^2.
\end{align}

\begin{remark}\label{cellassmp}
The inequality~\eqref{deltah2} implies that the distance between $\partial\Omega^n$ and $\partial\Omega^{n-1}$ is bounded by ${\cal O}(h^2)$. This implies the following property, which will be needed in the analysis below:
For each cell $K \in S_h^{n,n-1}$, there exists a path of cells $K_{i},\,i=1,\ldots,M$, such that $\overline{K}_{i} \cap \overline{K}_{i+1}$ is a facet in ${\cal F}_{h,\delta}^{g,n} \cap {\cal F}_{h,\delta}^{g,n-1}$ and the final cell $K_{M}$ lies fully in the interior of $\Omega^n$. Furthermore, the number of cases, in which an element $K_{M} \subset \Omega^{n}$ is utilised as a final element among all paths, can be bounded independently of $h$ and $\Delta t$.
\end{remark}

\subsubsection{Discrete variational formulation}
In the numerical approximation,  the boundary condition of the discrete problem (\ref{2FE1}) is implemented weakly by means of  Nitsche's method. Moreover, the function $u^{n}_{h}$ is extended by means of a ghost penalty term $g^{n}_{h}(\cdot,\cdot)$.
In each time step $n=1,2,\ldots,N$, we consider the following discrete variational formulation: 
Find $u^{n}_{h} \in V^{n,m}_{h}$ such that
\begin{align}\label{2FE1}
{\cal A}(u_h^n,u_h^{n-1}; v_h)
=(f^{n-\frac{1}{2}}_{h},v_{h}), \quad \forall v_h\in V^{n,m}_{h},
\end{align}
where 
\begin{align}\label{DefA}
{\cal A}(u_h^n,u_h^{n-1}; v_h) := \left(D_{\Delta t}^{-}u^{n}_{h},v_{h}
\right)_{\Omega^{n}}+\dfrac{1}{2}a^{n}_{h}(u^{n}_{h},v_{h})+\dfrac{1}{2}a^{n}_{h}(u^{n-1}_{h},v_{h})+\dfrac{\gamma_{D}}{h}(u^{n}_{h},v_{h})_{\partial \Omega^{n}}&+\gamma_{g}g^{n}_{h}(u^{n}_{h},v_{h})
\end{align}
and
\begin{align*}
D_{\Delta t}^{-}u^{n}_{h}=\dfrac{u^{n}_{h}-u^{n-1}_{h}}{\Delta t}, \quad 
a^{n}_{h}(u^{k}_{h},v_{h})=(\nabla u^{k}_{h},\nabla v_{h})_{\Omega^{n}}-(\partial_{n}u^{k}_{h},v_{h})_{\partial \Omega^{n}},\quad
f^{n-\frac{1}{2}}=\dfrac{f^{n}+f^{n-1}}{2}.
\end{align*}
We assume that all integrals in~\eqref{DefA} are evaluated exactly. For a consideration of additional quadrature errors, that result when cut cells are approximated linearly in the computation of the integrals, we refer to~\cite{Lc_esiam19}.

The ghost penalty stabilization is defined by 
\begin{equation}\label{ghostdef}
    g^{n}_{h}(w_{h},v_{h}) = \sum_{e \in \mathcal{F}^{n,g}_{h,\delta}} \sum_{k=1}^{m}\dfrac{h^{2k-1}}{k!^{2}} \int_e [[\partial^{k}_{n}w_{h}]]\cdot[[\partial^{k}_{n}v_{h}]]\,\text{d}s,
\end{equation}
where $[[\cdot]]$ is the jump operator and $\partial_{n}$ the exterior normal derivative. For further possibilities for the extension $g_h^n$, we refer to~\cite{Lc_esiam19}. The variant chosen here based on the jump of derivatives over edges has the advantage that it is fully consistent, in the sense that $g^{n}_{h}(u,v)$ vanishes for $u \in H^{m+1}(\Omega^{n}_{\delta})$.
The purpose of the ghost penalty is twofold: First, it serves to extend the solution $u_h^n$ implicitly to $\Omega_{h,\delta}^n$. Secondly, it ensures the discrete coercivity of the formulation~\eqref{2FE1} on ${\cal T}_{h,\delta}^n$. 

To incorporate the initial condition, we set $u_h^0 := E^1 u^0$ in~\eqref{DefA} for $n=1$, where $E^1 u^0$ is a smooth, e.g.~a  canonical extension, of the initial value $u_0$. This corresponds to the following Ritz projection of the initial value $u^0$
\begin{align*}
\Delta t^{-1} (u_h^0, v_h)_{\Omega^1} + a^1(u_h^0, v_h) = \Delta t^{-1} (E^1 u^0, v_h)_{\Omega^1} + a^1(E^1 u^0, v_h) \quad\forall v_h \in V_h^1,
\end{align*}

\noindent The following lemma is the key to extend the discrete coercivity to ${\Omega_\delta^n}$:
\begin{lemma}\label{extnbdd}
Given Assumption~\ref{ass:CFL}, any discrete function $v^{n}_{h}\in V^{n,m}_{h}$ satisfies 
\begin{align*}
    \|v_h^n\|^{2}_{\Omega_{h,\delta}^n} \leq c \|v_h^n\|^{2}_{\Omega^n}+ch^{2} g^{n}_{h}(v_{h},v_{h}),\qquad
    \|\nabla v_h^n\|^{2}_{\Omega_{h,\delta}^n} \leq c \|\nabla v_h^n\|^{2}_{\Omega^n}+c g^{n}_{h}(v_{h},v_{h}).
\end{align*}
In addition, for $v,w \in H^{m+1}(\Omega_\delta^n),\,m\geq 1$, it holds
\begin{align*}
    g^{n}_{h}(v,w) \leq g^{n}_{h}(v,v)^{1/2}g^{n}_{h}(w,w)^{1/2},\quad g^{n}_{h}(v,v) \leq c h^{2m}\| v\|^{2}_{H^{m+1}(\Omega_\delta^n)}.
\end{align*}
\end{lemma}
{\bf Proof.} A proof of this lemma is given in \cite{Lc_esiam19}.
\hfill \eop

\bigskip

At the end of this section, we briefly show that 
the variational formulation~\eqref{2FE1} is well-posed for each $n$. We define the discrete energy as
\begin{equation}\label{triplnorm}
    \mathfrak{E}^n(u^{n}_{h}, u^{n-1}_{h})=\left(\frac{1}{2}\|\nabla u^{n}_{h}+\nabla u^{n-1}_{h}\|_{\Omega^{n}}^{2}+\dfrac{1}{\Delta t}\|u^{n}_{h}-u^{n-1}_{h}\|_{\Omega^{n}}^{2}+ \dfrac{\gamma_{D}}{h}\|u^{n}_{h}\|^{2}_{\partial \Omega^{n}}+
\gamma_{g}g^{n}_{h}(u^{n}_{h},u^{n}_{h}) \right)^{1/2}.
\end{equation}
and the energy norm as
\begin{align*}
    |||u_h^n|||_n := \mathfrak{E}^n(u^n_h,0).
\end{align*}
We will show the coercivity relation
\begin{align}\label{coerc}
    |||u_h^n|||_n^2 \leq \dfrac{1}{4} {\cal A}(u_h^n,0; u_h^n), \qquad u_h^n\in V_h^{m,n}.
\end{align}
for sufficiently large $\gamma_g, \gamma_D$.
The well-posedness of~\eqref{2FE1} follows then by standard arguments.

From the definition of the bilinear form ${\cal A}(u_h^n,0; u_h^n) $, we have 
\begin{align*}
 {\cal A}(u_h^n,0; u_h^n)= \dfrac{1}{\Delta t}\|u^{n}_{h}\|_{\Omega^{n}}^2+\dfrac{1}{2}\|\nabla u^{n}_{h}\|_{\Omega^{n}}^2-\dfrac{1}{2}(\partial_{n}u^{n}_{h},u^{n}_{h})_{\partial \Omega^{n}}+\dfrac{\gamma_{D}}{h}\|u^{n}_{h}\|_{\partial \Omega^{n}}^2+\gamma_{g}g^{n}_{h}(u^{n}_{h},u^{n}_{h}). 
\end{align*}
The term $-\frac{1}{2}(\partial_{n}u^{n}_{h},u^{n}_{h})_{\partial \Omega^{n}}$ is estimated by means of Young's inequality, an inverse inequality and Lemma \ref{extnbdd} as follows for sufficiently small $\epsilon>0$:
\begin{align*}
    -\dfrac{1}{2}(\partial_{n}u^{n}_{h},u^{n}_{h})_{\partial \Omega^{n}} & \geq   -\dfrac{\epsilon h}{4 }\|\nabla u^{n}_{h}\|^{2}_{\partial \Omega^{n}}-\dfrac{1}{4 \epsilon h}\|u^{n}_{h}\|^{2}_{\partial \Omega^{n}}\\[4pt]
    & \geq   -c \epsilon \|\nabla u^{n}_{h}\|^{2}_{ \Omega^{n}_h}-\dfrac{1}{4 \epsilon h}\|u^{n}_{h}\|^{2}_{\partial \Omega^{n}}.
    \end{align*}
    For sufficiently large parameters $\gamma_D, \gamma_g$, we have
    \begin{align*}
    -\dfrac{1}{2}(\partial_{n}u^{n}_{h},u^{n}_{h})_{\partial \Omega^{n}}& \geq -\dfrac{1}{4} \left(\|\nabla u^{n}_{h}\|^{2}_{ \Omega^{n}}+\gamma_{g}g^{n}_{h}(u^{n}_{h},u^{n}_{h})\right)-\dfrac{\gamma_{D}}{4h}\|u^{n}_{h}\|^{2}_{\partial \Omega^{n}}.
\end{align*}
This proves the coercivity~\eqref{coerc}. 

\section{Stability analysis}\label{Sec:3}
In this section, a detailed stability analysis of the discrete problem (\ref{2FE1}) is developed. One of the main issues in the analysis is that the  discrete functions $u^{n-1}_{h}$ and $\nabla u^{n-1}_{h}$ appear on $\Omega^{n}$ in the $n$-th time step, whereas bounds are only available for $\| u^{n-1}_{h}\|_{\Omega^{n-1}}$ and  $\|\nabla u^{n-1}_{h}\|_{\Omega^{n-1}}$ from the previous time step. We start with some technical lemmas that will enable us to deal with this issue. 

\subsection{Auxiliary  estimates}

\begin{lemma}\label{tqlemma1}
Let Assumptions~\ref{ass:CFL} be valid.
Any discrete functions $v^{n}_{h}\in V^{n,m}_{h},\, v_h^{n-1}\in V_h^{n-1,m}$ satisfy the following inequality:
\begin{equation}\label{tql1}
    \|\nabla v^{n}_{h} - \nabla v^{n-1}_{h}\|^{2}_{\Omega^{n} \setminus \Omega^{n-1}} \leq \dfrac{2}{\Delta t}\| v^{n}_{h}- v^{n-1}_{h}\|^{2}_{\Omega^{n}}+ cg^{n}_{h}(v^{n}_{h},v^{n}_{h})+c g^{n-1}_{h}(v^{n-1}_{h},v^{n-1}_{h}).
\end{equation}
\end{lemma}
{\bf Proof.}
For $m=1$,$\nabla v_h^{n}-\nabla v_h^{n-1}$ is constant on a cell $T\in {\cal T}_h^n$. Due to the quasi-uniformity of the background mesh, we have $|T| \geq ch^d$ and from~\eqref{DeltaUpperBound} $|T \cap (\Omega^n\setminus \Omega^{n-1})| \leq ch^{d-1} \Delta t$. It follows that
\begin{align}\label{L31m1}
\|\nabla v^{n}_{h} - \nabla v^{n-1}_{h}\|^{2}_{\Omega^{n} \setminus \Omega^{n-1}} \leq c \dfrac{\Delta t}{h} \|\nabla v^{n}_{h} - \nabla v^{n-1}_{h}\|^{2}_{\mathcal{S}^{n,n-1}_{h} }.
\end{align}
In the general case ($m >1$), we have still 
\begin{align}\label{L31m2}
\|\nabla v^{n}_{h} - \nabla v^{n-1}_{h}\|^{2}_{\Omega^{n} \setminus \Omega^{n-1}} \leq \|\nabla v^{n}_{h} - \nabla v^{n-1}_{h}\|^{2}_{\mathcal{S}^{n,n-1}_{h} }.
\end{align}
Now let $K_{1} \in  \mathcal{S}^{n,n-1}_{h}$. By Remark~\ref{cellassmp},
there is a set of neighbouring cells $K_2,...,K_M$, such that
 $(\overline{K}_{i} \cap \overline{K}_{i+1}) \in ({\cal F}_{h,\delta}^{g,n} \cap {\cal F}_{h,\delta}^{g,n-1})$ and
 $K_M$ lies fully in the interior of $\Omega^n$. Let $e_{1,2}$ be the edge that separates the cells $K_{1}$ and $K_{2}$, see Figure~\ref{fig:domaindisp}. Then, using arguments from~\cite[Lemma 5.1]{massing2014stabilized} and~\cite[Lemma 5.2]{Lc_esiam19}, we can deduce that
\begin{equation}\label{l1tc1}
\begin{array}{ll}
 \|\nabla v^{n}_{h} - \nabla v^{n-1}_{h}\|^{2}_{K_1} &
 \displaystyle \leq c  \|\nabla v^{n}_{h} - \nabla v^{n-1}_{h}\|^{2}_{K_2}+c \sum_{k=1}^m \int_{e_{1,2}}  h^{2k-1} [[\nabla^k v^{n}_{h} - \nabla^k v^{n-1}_{h}]]^{2}ds\\[10pt]
 &
 \displaystyle \leq c  \|\nabla v^{n}_{h} - \nabla v^{n-1}_{h}\|^{2}_{K_2}+c \sum_{k=1}^m \int_{e_{1,2}}h^{2k-1} [[\nabla^k v^{n}_{h} ]]^{2}ds+c \int_{e_{1,2}}h^{2k-1}[[\nabla^k v^{n-1}_{h} ]]^{2}ds. 
 \end{array}
\end{equation}
We follow this process from $K_2$ to $K_M$ by crossing edges $e_{2,3}$ to $e_{M-1,M}$ to obtain 
\begin{align}
\begin{split}
\label{l1tc11}
    \|\nabla v^{n}_{h} - \nabla v^{n-1}_{h}\|_{K_1}^2 \leq &c   \|\nabla v^{n}_{h} - \nabla v^{n-1}_{h}\|^{2}_{K_M} \\
    &+ \sum_{k=1}^m \sum_{j=2}^M\left(\int_{e_{j-1,j}}h^{2k-1} [[\nabla^k v^{n}_{h} ]]^{2}ds+c \int_{e_{j-1,j}}h^{2k-1}[[\nabla^k v^{n-1}_{h} ]]^{2}ds\right).
    \end{split} 
\end{align}
For the first term on the right-hand side of~\eqref{l1tc11}, we use an inverse inequality and the CFL condition with a sufficiently small constant $c_{CFL}$
\begin{equation}\label{L31inveq}
\|\nabla v^{n}_{h} - \nabla v^{n-1}_{h}\|^{2}_{K_M}
 \leq \frac{c}{h^2}   \|v^{n}_{h} - v^{n-1}_{h}\|^{2}_{K_M} 
  \leq \frac{2}{\Delta t}   \|v^{n}_{h} - v^{n-1}_{h}\|^{2}_{K_M}.
\end{equation}
  As all edges $e_{j-1,j}$ belong to both ${\cal F}_{h,\delta}^{g,n}$ and ${\cal F}_{h,\delta}^{g,n-1}$, we obtain after summation over all cells in $\mathcal{S}^{n,n-1}_{h}$ from~\eqref{L31m2},~\eqref{l1tc11} and~\eqref{L31inveq}
  \begin{align*}
      \|\nabla v^{n}_{h} - \nabla v^{n-1}_{h}\|^{2}_{\Omega^{n} \setminus \Omega^{n-1}} 
      &\leq \frac{2}{\Delta t}   \|v^{n}_{h} - v^{n-1}_{h}\|^{2}_{\Omega^n} + c \left( g^{n}_{h}(v^{n}_{h},v^{n}_{h}) + g^{n-1}_{h}(v^{n-1}_{h},v^{n-1}_{h})\right).
  \end{align*}
  Using~\eqref{L31m1} instead of~\eqref{L31m2} for $m=1$, 
  the same result follows under the CFL condition $\Delta t \leq c_{CFL}h^{2}$ for an arbitrary constant $c_{CFL}$.
\hfill \eop
\bigskip

We note that the CFL condition is required to estimate the ``mismatch'' $\|\nabla v_h^{n-1}\|_{\Omega^n\setminus\Omega^{n-1}}^2$ by means of the discrete time derivative $\frac{1}{\Delta t}\|v_h^n-v_h^{n-1}\|_{\Omega^n}^2$ (see the following lemma).

\vskip 2mm

Next, we discuss how the term $\|\nabla u^{n-1}_{h}\|_{\Omega^{n}}$ can be bounded by $\|\nabla u^{n-1}_{h}\|_{\Omega^{n-1}}$ plus further terms that can be controlled in the following stability analysis.

\begin{lemma}\label{tqlemma2}
Under the assumptions of Lemma~\ref{tqlemma1} it holds for $v_h^n \in V_h^{n,m}, \,v_h^{n-1}\in V_h^{n-1,m}$ that
\begin{equation}\label{tql2}
\begin{array}{ll}
\Delta t\|\nabla v^{n-1}_{h}\|^{2}_{\Omega^{n}\setminus\Omega^{n-1}}& \leq 
 \dfrac{\Delta t}{2} \|\nabla v^{n}_{h}+\nabla v^{n-1}_{h}\|^{2}_{\Omega^{n}\setminus \Omega^{n-1}}+ \| v^{n}_{h}- v^{n-1}_{h}\|^{2}_{\Omega^{n}}\\[10pt]
&\qquad + c \Delta t g^{n}_{h}(v^{n}_{h},v^{n}_{h})+c  \Delta t g^{n-1}_{h}(v^{n-1}_{h},v^{n-1}_{h}).
\end{array}
\end{equation}
\end{lemma}
{\bf Proof.} 
By the triangle inequality, we have 
\begin{align*}
\|\nabla v^{n-1}_{h}\|_{\Omega^{n} \setminus \Omega^{n-1}} \leq \dfrac{1}{2}\|\nabla v^{n}_{h}+\nabla v^{n-1}_{h}\|_{\Omega^{n} \setminus \Omega^{n-1}}+\dfrac{1}{2}\|\nabla v^{n}_{h}-\nabla v^{n-1}_{h}\|_{\Omega^{n} \setminus \Omega^{n-1}}.
\end{align*}
By means of the inequality $(a+b)^2 \leq 2a^2 + 2b^2$ this implies
\begin{equation}\label{l2tc2}
\Delta t\|\nabla v^{n-1}_{h}\|^{2}_{\Omega^{n} \setminus \Omega^{n-1}} \leq \dfrac{\Delta t}{2}\|\nabla v^{n}_{h}+\nabla v^{n-1}_{h}\|^{2}_{\Omega^{n} \setminus \Omega^{n-1}}+\dfrac{\Delta t}{2}\|\nabla v^{n}_{h}-\nabla v^{n-1}_{h}\|^{2}_{\Omega^{n} \setminus \Omega^{n-1}}.
\end{equation}
The statement follows by using Lemma~\ref{tqlemma1} for the second term in~\eqref{l2tc2}
\begin{align}\label{l2tc4}
\begin{split}
\dfrac{\Delta t}{2}\|\nabla v^{n}_{h}-\nabla v^{n-1}_{h}\|^{2}_{\Omega^{n} \setminus \Omega^{n-1}} &\leq   \|v^{n}_{h} - v^{n-1}_{h}\|^{2}_{\Omega^{n}} + c \Delta t g^{n}_{h}(v^{n}_{h},v^{n}_{h})+c \Delta t g^{n-1}_{h}(v^{n-1}_{h},v^{n-1}_{h}). 
\end{split}
\end{align}
\hfill \eop

\noindent Next, we provide the following Poincar\'e-type estimate:
\begin{lemma}\label{lem:poincare}
Let $u\in W^{1,p}(\Omega)$ for $1\leq p<\infty$ and let the CFL condition (Assumption~\ref{ass:CFL}) be valid. It holds for $l\in \{n-1,n\}$
\begin{align}
\|u\|_{L^p(\Omega^n \setminus \Omega^{n-1})}^p \leq c\Delta t \|u\|_{L^p(\partial\Omega^l)}^p + c\Delta t^2 \|\nabla u\|_{L^p(\Omega^n \setminus \Omega^{n-1})}^p.
\end{align}
\end{lemma}
{\bf Proof.} The proof follows the lines of~\cite[Lemma 4.34]{richter2017fluid} and uses the fact that dist$(\Omega^n, \Omega^{n-1}) \leq c\Delta t$.
\eop

\medskip
\noindent Using this, we can derive bounds for $\| v^{n-1}_{h}\|^{2}_{\Omega^{n}\setminus\Omega^{n-1}}$ and $\|v^{n}_{h}\|^{2}_{\Omega^{n} \setminus \Omega^{n-1}}$:
\begin{lemma}\label{tqlemma3}
Under the assumptions of Lemma~\ref{tqlemma1} it holds for $v^{l}_{h}\in V^{l,m}_{h}, \, l \in \{n-1, n\}$
\begin{eqnarray}
 \| v^{l}_{h}\|^{2}_{\Omega^{n}\setminus\Omega^{n-1}} &\leq c \Big( \Delta t 
 \| v^{l}_{h}\|^{2}_{\partial \Omega^{l}}
 + \Delta t \| v^{n}_{h}- v^{n-1}_{h}\|^{2}_{\Omega^{n}} 
+ \Delta t^2 \|\nabla v^{n}_{h}+\nabla v^{n-1}_{h}\|^{2}_{\Omega^{n}}\nonumber \\
 &\hspace{2cm}+\Delta t^{2}g^{n-1}_{h}( v^{n-1}_{h}, v^{n-1}_{h})
+ \Delta t^{2}g^{n}_{h}( v^{n}_{h}, v^{n}_{h})\Big) 
 \label{tql3}
\end{eqnarray}

\end{lemma}
{\bf Proof.}
By means of Lemma~\ref{lem:poincare} for $p=2$ we have for $l\in \{n-1,n\}$
\begin{align}\label{l3tc1}
\| v^{l}_{h}\|^{2}_{\Omega^{n} \setminus \Omega^{n-1}}& \leq c \Delta t \| v^{l}_{h}\|^{2}_{\partial \Omega^{l} }+c\Delta t^{2}\|\nabla  v^{l}_{h}\|^{2}_{\Omega^{n} \setminus \Omega^{n-1}}
\end{align}
The statement follows by applying Lemma~\ref{tqlemma2} to the second term in~\eqref{l3tc1}.
\hfill \eop

\subsection{Stability result}

Before discussing the stability result, we introduce some abbreviations for the space-time Bochner norms to  simplify the  mathematical expressions 
\[
\|u\|_{\infty,m,I_k}:=\|u\|_{L^{\infty}(I_{k},H^{m}(\Omega(t)))}, \qquad \|u\|_{\infty,m}:=\|u\|_{\infty,m,I},
\]
where $m \in \mathbb{N}\cup\{0\}$ and $H^{0}(\Omega(t)):=L^{2}(\Omega(t))$.

Now we are ready to prove the following stability result.

\begin{theorem}[Stability]\label{theo.stab}
Let Assumptions~\ref{ass:CFL} be valid and let $f\in L^{\infty}(I, L^2(\Omega(t))),\, u^0\in H^1(\Omega^0)$ and let the mapping $T$ be a $W^{1,\infty}$-diffeomorphism. For sufficiently large $\gamma_g, \gamma_D$ the solution $\{u^{k}_{h}\}_{k=1}^n$ of the discrete problem (\ref{2FE1}) fulfills
   \begin{equation}\label{stbeq}
\|u^{n}_{h}\|^{2}_{\Omega^{n}}+\Delta t \|u^{n}_{h}\|^{2}_{\Omega^{n}} +\Delta t \sum_{k=1}^{n} \mathfrak{E}^k(
u^{k}_{h}, u^{k-1}_{h})
\leq c \exp(c \,t_{n}) \left( \|u^{0}\|^{2}_{\Omega^{0}}+t_{n}\|\nabla u^{0}\|^{2}_{\Omega^{0}}+t_{n} \|f\|^{2}_{\infty,0}\right).
\end{equation} 
\end{theorem}
{\bf Proof.}
We test (\ref{2FE1}) with $v_{h}=2 \Delta t u^{n}_{h}$ to obtain 
\begin{equation}\label{2stb1}
\begin{array}{ll}
2\left(u^{n}_{h}-u^{n-1}_{h}, u^{n}_{h}\right)_{\Omega^{n}}+\Delta t \left(\nabla u^{n}_{h}+\nabla u^{n-1}_{h}, \nabla u^{n}_{h}\right)_{\Omega^{n}} +\dfrac{2\gamma_{D} \Delta t}{ h}(u^{n}_{h},u^{n}_{h})_{\partial \Omega^{n}}\\[16pt]
- \Delta t (\partial_{n}u^{n}_{h}+\partial_{n}u^{n-1}_{h},u^{n}_{h})_{\partial \Omega^{n}}+2\Delta t\gamma_{g}g^{n}_{h}(u^{n}_{h},u^{n}_{h})=2\Delta t(f^{n-\frac{1}{2}},u^{n}_{h}).
\end{array}
\end{equation}
We estimate the fourth term in (\ref{2stb1}) by means of Young's inequality with a sufficiently small $\epsilon>0$ followed by an inverse inequality
\begin{align*}
\Delta t (\partial_{n}u^{n}_{h}+\partial_{n}u^{n-1}_{h},u^{n}_{h})_{\partial \Omega^{n}}
&\leq  \frac{\Delta t \epsilon h}{16} \| \partial_n u^{n}_{h}+\partial_n u^{n-1}_{h}\|^{2}_{\partial \Omega^n}+\dfrac{4\Delta t}{\epsilon h}  \|u^{n}_{h}\|_{\partial\Omega^n}^2\\[8pt]
&\leq  \dfrac{\Delta t }{16} \| \nabla u^{n}_{h}+\nabla u^{n-1}_{h}\|^{2}_{ \Omega^n}+\dfrac{\gamma_D \Delta t}{2 h}  \|u^{n}_{h}\|_{\partial\Omega^n}^{2},
\end{align*}
 where $\gamma_D \geq 8/\epsilon$. By using the relation $2(a+b,a)=(a+b)^{2}+a^{2}-b^{2}$ for the first two terms in~\eqref{2stb1}, we obtain that 
\begin{align}\label{2stb2}
\begin{split}
\|u^{n}_{h}\|^{2}_{\Omega^{n}}+\|u^{n}_{h}-&u^{n-1}_{h}\|^{2}_{\Omega^{n}}+\dfrac{\Delta t}{2}\|\nabla u^{n}_{h}\|_{\Omega^{n}}^{2}+\dfrac{7\Delta t}{16}\|\nabla u^{n}_{h}+\nabla u^{n-1}_{h}\|_{\Omega^{n}}^{2}
+\dfrac{3\gamma_{D} \Delta t}{2 h}\|u^{n}_{h}\|^{2}_{\partial \Omega^{n}} \\
&+2\Delta t\gamma_{g}g^{n}_{h}(u^{n}_{h},u^{n}_{h}) \,\leq\, \|u^{n-1}_{h}\|^{2}_{\Omega^{n}}+\dfrac{\Delta t}{2}\|\nabla u^{n-1}_{h}\|^{2}_{\Omega^{n}}+2\Delta t(f^{n-1/2},u^{n}_{h}).
\end{split}
\end{align}
For $n>1$, we bring the terms $\|u_h^{n-1}\|_{\Omega^n}$ and $\|\nabla u_h^{n-1}\|_{\Omega^n}$ to the domain $\Omega^{n-1}$.
By employing Lemmas~\ref{tqlemma2} and~\ref{tqlemma3},  we have
\begin{align}\label{2stb3}
\begin{split}
 \|u^{n-1}_{h}\|^{2}_{\Omega^{n}}+\dfrac{\Delta t}{2}\|\nabla u^{n-1}_{h}\|^{2}_{\Omega^{n}}
 &\leq \| u^{n-1}_{h}\|^{2}_{\Omega^{n-1}}+c \Delta t\| u^{n-1}_{h}\|^{2}_{\partial \Omega^{n-1}}
 +\dfrac{\Delta t}{2}\|\nabla u^{n-1}_{h}\|^{2}_{\Omega^{n-1}}\\
 &\quad+ \left(\frac{\Delta t}{4}+c\Delta t^2\right) \|\nabla u^{n}_{h}+\nabla u^{n-1}_{h}\|^{2}_{\Omega^{n}}+ \left(\frac{1}{2} +c\Delta t\right) \|u^{n}_{h}- u^{n-1}_{h}\|^{2}_{\Omega^{n}}\\
&\qquad+ c \Delta t g^{n}_{h}(u^{n}_{h},u^{n}_{h})+c  \Delta t g^{n-1}_{h}(u^{n-1}_{h},u^{n-1}_{h}).
\end{split}
\end{align}
Inserting (\ref{2stb3}) into (\ref{2stb2}) and using $2\Delta t (f^{n-1/2}, u_h^n)_{\Omega^n} \leq \Delta t \left(\|f^{n-1/2}\|_{\Omega^n}^2 + \|u_h^n\|_{\Omega^n}^2\right)$ gives for sufficiently large $\gamma_{g}$
\begin{align}\label{2stb5}
    \begin{split}
(1-\Delta t)\|&u^{n}_{h}\|^{2}_{\Omega^{n}}+\dfrac{\Delta t}{2}\|\nabla u^{n}_{h}\|_{\Omega^{n}}^{2}+ \dfrac{1}{4}\|u^{n}_{h}-u^{n-1}_{h}\|^{2}_{\Omega^{n}}+\dfrac{\Delta t}{8}\|\nabla u^{n}_{h}+\nabla u^{n-1}_{h}\|_{\Omega^{n}}^{2}
\\
&+\dfrac{3\gamma_{D} \Delta t}{2 h} \|u^{n}_{h}\|^{2}_{\partial \Omega^{n}} + \Delta t\gamma_{g}g^{n}_{h}(u^{n}_{h},u^{n}_{h})\\
&\qquad\qquad \leq \| u^{n-1}_{h}\|^{2}_{\Omega^{n-1}}+\dfrac{\Delta t}{2}\|\nabla u^{n-1}_{h}\|^{2}_{\Omega^{n-1}} +c \Delta t\| u^{n-1}_{h}\|^{2}_{\partial \Omega^{n-1}}\\
&\qquad\qquad\qquad\qquad+c \Delta t g^{n-1}_{h}( u^{n-1}_{h}, u^{n-1}_{h})+ \Delta t \|f^{n-\frac{1}{2}}_{h}\|^{2}_{\Omega^{n}}.
\end{split}
\end{align}
For $n=1$, we obtain from (\ref{2stb2}) and the stability of the extension $E^1$
\begin{align}
    \begin{split}
(1-\Delta t)\|&u^{1}_{h}\|^{2}_{\Omega^{1}}+ \|u^{1}_{h}-u^{0}_{h}\|^{2}_{\Omega^{1}}+\dfrac{\Delta t}{2}\|\nabla u^{1}_{h}\|_{\Omega^{1}}^{2} + \dfrac{3\Delta t}{8}\|\nabla u^{1}_{h}+\nabla u^{0}_{h}\|_{\Omega^{1}}^{2}
\\
&+\dfrac{3\gamma_{D} \Delta t}{2 h} \|u^{1}_{h}\|^{2}_{\partial \Omega^{1}} + 2\Delta t\gamma_{g}g^{1}_{h}(u^{1}_{h},u^{1}_{h})\\
&\qquad\qquad \leq \| E^1 u^0\|^{2}_{\Omega^{1}}+\dfrac{\Delta t}{2}\|\nabla (E^1u^{0})\|^{2}_{\Omega^{1}} + \Delta t \|f^{\frac{1}{2}}_{h}\|^{2}_{\Omega^{1}}\\
&\qquad\qquad\qquad\qquad \leq c\| u^0\|^{2}_{\Omega^{0}}+c \Delta t\|\nabla u^{0}\|^{2}_{\Omega^{0}} + \Delta t \|f^{\frac{1}{2}}_{h}\|^{2}_{\Omega^{1}}.
\end{split}
\end{align}
Taking the sum over $k=1,2,\ldots,n$, this yields for sufficiently large $\gamma_g$ and $\gamma_D$ 
\begin{align}\label{2stb6}
\begin{split}
\displaystyle \|u^{n}_{h}\|^{2}_{\Omega^{n}}+\Delta t\|\nabla u^{n}_{h}\|_{\Omega^{n}}^{2}&+\dfrac{\Delta t}{4}\sum_{k=1}^{n}
\mathfrak{E}^k(u^{k}_{h},u^{k-1}_{h}) \\
&\leq \displaystyle c \|u^{0}_{h}\|^{2}_{\Omega^{0}}+c t_{n}\|\nabla u^{0}_{h}\|^{2}_{\Omega^{1}}+2t_{n} \|f\|^{2}_{\infty,0}+ c \Delta t \sum_{k=2}^{n} \| u^{k-1}_{h}\|^{2}_{\Omega^{k-1}}.
\end{split}
\end{align}
The statement follows by means of the discrete Gronwall lemma.
\hfill \eop

\begin{remark}[CFL condition]
For $m=1$, the stability result in Theorem~\ref{theo.stab} could be shown under the weaker CFL condition $\Delta t \leq c_{CFL} h^{3/2}$ for sufficiently small $c_{CFL}$. This is due to the possibility to use the estimate
\begin{align}
\|\nabla  v^{n-1}_{h}\|^{2}_{\Omega^{n} \setminus \Omega^{n-1}} \leq
c\frac{\Delta t}{h}\|\nabla  v^{n-1}_{h}\|^{2}_{ \Omega^{n-1}_{h,\delta}}
\end{align}
in~\eqref{L31m1} (Note that, for $m=1$, $\nabla  v^{n-1}_{h}$ is constant in each cell $T$). However, in the following section the stronger CFL condition $\Delta t \leq c_{CFL} h^2$ with arbitrary $c_{CFL}$ (see Assumption~\ref{ass:CFL}) will be used in order to show optimal-order convergence estimates. The CFL condition is needed to estimate the term $\|\nabla v_h^{n-1}\|_{\Omega^n \setminus \Omega^{n-1}}^2$ by the discrete time derivative $\Delta t^{-1} \|v_h^n - v_h^{n-1}\|_{\Omega^n}^2$, see~\eqref{L31inveq} and Lemma~\ref{tqlemma2}. 
\end{remark}

\section{A priori error analysis}\label{Sec:4}
In this section, we show an \textit{a priori} error estimate for the discrete problem (\ref{2FE1}). We define the discretisation error as 
\begin{equation}\label{defglberr}
 e^{n}:=u^{n}-u^{n}_{h}, \quad n \geq 1,
\end{equation}
where $u^{n}:=u(t_n)$ is assumed to be at least in $H^2(\Omega^n)$. Further regularity assumptions on $u$ will be made below.
The error  is decomposed into an interpolation error $\eta^{n}$ and a discrete error $\xi^{n}_{h}$ terms defined by
\begin{equation}\label{deferr}
\eta^{n}:=u^{n}-I^{n}_{h}u^{n},\qquad \xi^{n}_{h}:=I^{n}_{h}u^{n}-u^{n}_{h},
\end{equation}
where $ I^{n}_{h}u^{n} $ denotes the standard Lagrangian nodal interpolation of $u^{n}$ on $\mathcal{T}^{n}_{h,\delta}$. For $n=0$ we have, by definition of $u_h^0$, that $e^0 = u^0 - u_h^0=0$ in $\Omega^0$ and thus, we also set $\eta^0=\xi_h^0=0$. 
We will make use of the following standard interpolation estimates for $n\geq 1$
\begin{align}
    \|\eta^{n}\|_{H^{l}(\Omega)} \leq c h^{k-l}\|u^{n}\|_{H^{k}(\Omega)},\quad \mbox{for} \,\,0 \leq l \leq 1,\,\,2 \leq k \leq m+1, \label{interr1}\\
     \|\eta^{n}\|_{H^{l}(\partial \Omega)} \leq c h^{k-l-1/2}\|u^{n}\|_{H^{k}(\Omega)},\quad \mbox{for} \,\,0 \leq l \leq 1,\,\, 2 \leq k \leq m+1. \label{interr2}
\end{align}


\subsection{Consistency error}

The exact solution $u \in H^{1}(\Omega(t))$ of the continuous problem (\ref{C1}) satisfies the following weak formulation:
\begin{equation}\label{Cweak1}
(u_{t},v)_{\Omega(t)}+a(u,v)=(f,v)_{\Omega(t)},\quad t \in I_{n}
\end{equation}
for $v \in H^{1}(\Omega(t))$ and the bilinear form
\[
a(u,v)=(\nabla u, \nabla v)_{\Omega(t)}-(\partial_{n}u, v)_{\partial \Omega(t)}.
\]
At time $t_{n-1}$, we have
\begin{equation}\label{Cweak1n-1}
(u_{t}(t_{n-1}),v)_{\Omega^{n-1}}+a^{n-1}(u^{n-1},v)=(f^{n-1},v)_{\Omega^{n-1}}
\end{equation}
where 
\[
a^{n-1}(u^{n-1},v)=(\nabla u^{n-1}, \nabla v)_{\Omega^{n-1}}-(\partial_{n}u^{n-1}, v)_{\partial \Omega^{n-1}}.
\]

To estimate the consistency error, we will need an analogous equality for $u^{n-1}$ on $\Omega^n$ instead of $\Omega^{n-1}$. Therefore, we define $\tilde{u}$ as a smooth extension of the exact solution $u$ to $Q^{n}_{\delta}$. Moreover, we define a smooth extension of the source term $f$ as follows:
\begin{equation}\label{fext}
\tilde{f}(t_{n-1})=\tilde{u}_{t}(t_{n-1})-\Delta \tilde{u}(t_{n-1}),\quad \mbox{on} \,\,\Omega^{n}_{\delta}.
\end{equation}
It holds 
\[
\tilde{f}(t_{n-1}) = f(t_{n-1})\quad \mbox{on} \,\,\Omega^{n-1}
\]
However, as we have allowed an arbitrary smooth extension of $f$ to $\Omega_\delta^{n-1}$ in Section~\ref{Sec:2}, this does
not necessarily hold in $\Omega_\delta^{n-1} \setminus \Omega^{n-1}$.
By using a test function $v_{h}\in V_h^{n,m}$, we get
\begin{equation}\label{extvar1}
  (\tilde{u}^{n-1}_{t},v_{h})_{\Omega^{n}}+(\nabla \tilde{u}^{n-1},\nabla v_{h})_{\Omega^{n}}-(\partial_{n}\tilde{u}^{n-1}_{h},v_{h})_{\partial \Omega^{n}}=(\tilde{f}^{n-1},v_{h})_{\Omega^{n}}.
\end{equation}
By adding the equations (\ref{Cweak1}) and (\ref{extvar1}), we obtain 
\begin{equation}\label{extvar2}
  (u^{n}_{t}+\tilde{u}^{n-1}_{t},v_{h})_{\Omega^{n}}+(\nabla u^{n}+\nabla \tilde{u}^{n-1}, \nabla v_{h})_{\Omega^{n}}-(\partial_{n}u^{n}_{h}+\partial_{n}\tilde{u}^{n-1}_{h},v_{h})_{\partial \Omega^{n}}=(f^{n}+\tilde{f}^{n-1},v_{h})_{\Omega^{n}}. 
\end{equation}
We note that the right-hand side in~\eqref{extvar2} differs from the discrete formulation~\eqref{2FE1} by $(\tilde{f}^{n-1}-f^{n-1},v_{h})_{\Omega^{n}\setminus \Omega^{n-1}}$. Hence, (\ref{extvar2}) can be rewritten as 
\begin{equation}\label{extvar4}
\left(\dfrac{u^{n}_{t}+\tilde{u}^{n-1}_{t}}{2},v_{h}\right)_{\Omega^{n}}+\dfrac{1}{2}a^{n}_{h}(u^{n},v_{h})+\dfrac{1}{2}a^{n-1}_{h}(\tilde{u}^{n-1},v_{h})=\dfrac{1}{2}(f^{n}+f^{n-1},v_{h})_{\Omega^{n}}+\dfrac{1}{2} \mathcal{E}^{n-1}_{f}(v_h),
\end{equation}
where $ \mathcal{E}^{n-1}_{f}(v_h) $ is given by
\begin{equation}\label{Ef}
\mathcal{E}^{n-1}_{f}(v_h) = (\tilde{f}^{n-1}-f^{n-1},v_{h})_{\Omega^{n} \setminus \Omega^{n-1}}.
\end{equation}
In the next lines, we will discuss a bound for the term $\mathcal{E}^{n-1}_{f}$.
\begin{lemma}\label{Tdomainbd}
Under the assumptions of Lemma~\ref{tqlemma1}, the error term $\mathcal{E}^{n-1}_{f}$ defined in~\eqref{Ef} satisfies the following estimate for $v^{n}_h\in V_h^{n,m}$ and sufficiently regular $u$:
\begin{align}\label{fbddeq}
    	\Delta t\big|\mathcal{E}^{n-1}_{f}(v^{n}_h)\big| &\leq  c \Delta t^{5}   {\cal R}_C(u)^{2}+\dfrac{1}{32} \Big(\Delta t  \|\nabla v^{n}_{h}+\nabla v^{n-1}_{h}\|^{2}_{\Omega^{n} }+\| v^{n}_{h}- v^{n-1}_{h}\|^{2}_{\Omega^{n}}
\\[12pt]
& \quad + \Delta t^{2}  \| v^{n}_{h}\|^{2}_{\partial \Omega^{n} }
+   \Delta t  g^{n}_{h}(v^{n}_{h},v^{n}_{h})+ \Delta t  g^{n-1}_{h}(v^{n-1}_{h},v^{n-1}_{h})\Big),\nonumber
\end{align}
where ${\cal R}_C(u) = \|u\|_{\infty, 3, I} + \|\partial_t u\|_{\infty, 1, I}.$
\end{lemma}
\noindent
{\bf Proof.} We apply Lemma~\ref{lem:poincare} to $u=(f^{n-1} - \tilde{f}^{n-1})v_h^n$ for $p=1$ 
\begin{equation}\label{fteq1}
\begin{array}{ll}
\Delta t\big|(f^{n-1}-\tilde{f}^{n-1},v^{n}_{h})_{\Omega^{n} \setminus \Omega^{n-1}}\big| &\leq \Delta t\|(f^{n-1}-\tilde{f}^{n-1})v^{n}_{h}\|_{L^{1}(\Omega^{n} \setminus \Omega^{n-1})}\\[8pt]
&\leq c \Delta t^{2} \|(f^{n-1}-\tilde{f}^{n-1})v^{n}_{h}\|_{L^{1}(\partial \Omega^{n-1}) }\\[8pt]
&\qquad+c \Delta t^{3} \left\|\nabla\big[(f^{n-1}-\tilde{f}^{n-1}) v^{n}_{h}\big]\right\|_{L^{1}(\Omega^{n} \setminus \Omega^{n-1})}.
\end{array}
\end{equation}
As $f^{n-1}=\tilde{f}^{n-1}$ on $\partial \Omega^{n-1}$, the first term in~\eqref{fteq1} vanishes. By definition of $\tilde{f}$ in (\ref{fext}) we can estimate further
\begin{equation}\label{fteq11}
\begin{array}{ll}
\Delta t\big|(f^{n-1}-\tilde{f}^{n-1},v^{n}_{h})_{\Omega^{n} \setminus \Omega^{n-1}}\big| &\leq c \Delta t^{3} \left[ \|(f^{n-1}-\tilde{f}^{n-1})\|_{\Omega^{n} \setminus \Omega^{n-1}}\|\nabla v^{n}_{h}\|_{\Omega^{n} \setminus \Omega^{n-1}}\right]
\\[8pt]
&\quad  +c \Delta t^{3} \left[ \|\nabla(f^{n-1}-\tilde{f}^{n-1})\|_{\Omega^{n} \setminus \Omega^{n-1}}\| v^{n}_{h}\|_{\Omega^{n} \setminus \Omega^{n-1}}\right] \\[8pt]
& \leq c \Delta t^{3}  {\cal R}_C(u)\left[\|\nabla v^{n}_{h}\|_{\Omega^{n} \setminus \Omega^{n-1}}+\| v^{n}_{h}\|_{\Omega^{n} \setminus \Omega^{n-1}}\right].
\end{array}
\end{equation}
Now, by employing Lemmas \ref{tqlemma2} and~\ref{tqlemma3}, we obtain
\begin{equation}\label{fteq2}
\begin{array}{ll}
\|\nabla v^{n}_{h}\|^{2}_{\Omega^{n} \setminus \Omega^{n-1}} + \|v^{n}_{h}\|^{2}_{\Omega^{n} \setminus \Omega^{n-1}} 
& \leq \left(\dfrac{1}{2}+c\Delta t^2\right) \|\nabla v^{n}_{h}+\nabla v^{n-1}_{h}\|^{2}_{\Omega^{n}\setminus \Omega^{n-1}}+ \dfrac{c}{\Delta t} \| v^{n}_{h}- v^{n-1}_{h}\|^{2}_{\Omega^{n}}\\
&\qquad+c\Delta t \|v_h^n\|_{\partial\Omega^n}^2 + c g^{n}_{h}(v^{n}_{h},v^{n}_{h})+c  g^{n-1}_{h}(v^{n-1}_{h},v^{n-1}_{h}).
\end{array}
\end{equation}
Inserting these estimates into~\eqref{fteq11} and using Young's inequality and Assumption~\ref{ass:CFL}, we obtain
\begin{equation}\label{fteq3}
\begin{array}{ll}
&\Delta t\big|(f^{n-1}-\tilde{f}^{n-1},v^{n}_{h})_{\Omega^{n} \setminus \Omega^{n-1}}\big| 
 \leq c  \Delta t^{3}  {\cal R}_C(u) \left[\|\nabla v^{n}_{h}\|_{\Omega^{n} \setminus \Omega^{n-1}}+\| v^{n}_{h}\|_{\Omega^{n} \setminus \Omega^{n-1}}\right]\\[12pt]
&\qquad\qquad \leq c \Delta t^{5}  {\cal R}_C(u)^2+\dfrac{1}{32}\Big(\Delta t \|\nabla v^{n}_{h}+\nabla v^{n-1}_{h}\|^{2}_{\Omega^{n} }+\| v^{n}_{h}- v^{n-1}_{h}\|^{2}_{\Omega^{n}}
+ \Delta t^{2}\| v^{n}_{h}\|^{2}_{ \partial \Omega^{n} }
\\[12pt]
& \hspace*{5.6cm}
+  \Delta t \gamma_g g^{n}_{h}(v^{n}_{h},v^{n}_{h})+ \Delta t \gamma_g g^{n-1}_{h}(v^{n-1}_{h},v^{n-1}_{h})\Big).
\end{array}
\end{equation}
\hfill \eop

\bigskip

Now we are ready to estimate the consistency error related to the discrete problem (\ref{2FE1}). 
By subtracting (\ref{2FE1}) from (\ref{extvar4}), the global  error term $e^{n}$ satisfies the equality
\begin{equation}\label{Err1}
\left(D^{-}_{\Delta t} e^{n},v_{h}
\right)_{\Omega^{n}}+\dfrac{1}{2}a^{n}_{h}(e^{n},v_{h})+\dfrac{1}{2}a^{n-1}_{h}(e^{n-1},v_{h})+\dfrac{\gamma_{D}}{h}(e^{n},v_{h})_{\partial \Omega^{n}}+\gamma_{g}g^{n}_{h}(e^{n},v_{h})
=\dfrac{1}{2}\mathcal{E}^{n-1}_{f}(v_h)+\mathcal{E}^{n}_{c}(v_h),
\end{equation}
where the consistency error $\mathcal{E}^{n}_{c}(v_h)$ is given by
\begin{equation}\label{conseq}
\begin{array}{ll}
\mathcal{E}^{n}_{c}(v_h)&=\underbrace{\left(D^{-}_{\Delta t}u^{n}-\partial_{t}E^{n}u(t_{n-1/2}),v_{h}\right)_{\Omega^{n}}}_{I_{1}}-\underbrace{\left(\dfrac{u_{t}(t_{n})+\partial_{t}E^{n}u(t_{n-1})}{2}-\partial_{t}E^{n} u(t_{n-1/2}),v_{h}\right)_{\Omega^{n}}}_{I_{2}}\\[12pt]
& \qquad +\underbrace{\dfrac{\gamma_{D}}{h}(E^n u^{n},v_{h})_{\partial \Omega^{n}}}_{I_{3}}+\underbrace{\gamma_{g}g^{n}_{h}(E^n u^{n},v_{h})}_{I_{4}}.
\end{array}
\end{equation}

The terms $I_3$ and $I_4$ vanish due to the homogeneous boundary condition and the regularity assumption on the exact solution $u^n\in H^2(\Omega^n)$ and its extension $E^n u^n\in H^2(\Omega_\delta^n)$.
The remaining terms are estimated in the following lemma.
\begin{lemma}\label{cnstlemma}
Let $u \in W^{3,\infty}(Q^n)$. Under Assumption~\ref{ass:regularity}, the consistency error for $v_{h} \in V^{n,m}_{h}$ is bounded by
\begin{equation}\label{consteq}
 \Delta t|\mathcal{E}^{n}_{c} (v^{n}_{h})| \leq   c   \Delta t^{5}  \|u\|^{2}_{W^{3,\infty}(Q^{n})}+\dfrac{\Delta t}{16}\|v^{n}_{h}\|^{2}_{\Omega^{n}}. 
\end{equation}
\end{lemma}
\noindent
{\bf Proof.}
First, we will show a bound for the term $I_{1}$.  By following the argumentation in \cite[Chapter~1]{thomee2007galerkin}, we have
\[
\begin{array}{ll}
\left|\left(\dfrac{u^{n}-E^{n}u^{n-1}}{\Delta t}-\partial_{t}E^{n}u(t_{n-1/2})\right)\right| &  \leq 
\displaystyle \dfrac{1}{\Delta t} \left|\int_{t_{n-1/2}}^{t^{n}}\dfrac{(t_{n}-s)^{2}}{3!}\partial^{3}_{t}E^{n}u(s)ds\right| \\ [16pt]
&  
\displaystyle \quad +\dfrac{1}{\Delta t} \left|\int^{t_{n-1/2}}_{t^{n-1}}\dfrac{(s-t_{n-1})^{2}}{3!}\partial^{3}_{t}E^{n}u(s)ds\right| \\ [16pt]
& 
\displaystyle \leq c \Delta t^{2} \sup_{ t \in [t_{n-1},t_{n}]}|\partial^{3}_{t}E^{n}u(t)|.
\end{array}
\]
Using the stability of the extension operator $E^{n}$ given in \eqref{exteq3} and the Cauchy-Schwarz inequality, we have
\begin{align}\label{cons1}
\begin{split}
\left|\left(\dfrac{u^{n}-E^{n}u^{n-1}}{\Delta t}-\partial_{t}E^{n}u(t_{n-1/2}),v^{n}_{h}\right)_{\Omega^{n}}\right| &\leq c \Delta t^{2} \|\partial^{3}_{t}E^{n}u^{n}\|_{L^{\infty}(Q^{n}_{\delta})}\|v^{n}_{h}\|_{\Omega^{n}} \\
&\leq  c \Delta t^{2} \|u\|_{W^{3,\infty}(Q^{n})}\|v^{n}_{h}\|_{\Omega^{n}}.
\end{split}
\end{align}
A bound for the second term $I_{2}$ follows in a similar way, see also~\cite[Chapter~1]{thomee2007galerkin}. The statement follows using Young's inequality.
\hfill \eop

\subsection{Interpolation error}
To derive an interpolation error estimate, we devise a discrete problem associated with the discrete error $\xi^{n}_{h}$. By definition of $\xi^{n}_{h}$  (\ref{deferr}) and  using (\ref{Err1}), we have for $v_h\in V_h^{n,m}$
\begin{equation}\label{Err2}
\begin{array}{ll}
\hspace{-4cm}\left(D^{-}_{\Delta t}\xi^{n}_{h},v_{h}
\right)_{\Omega^{n}}+\dfrac{1}{2}a^{n}_{h}(\xi^{n}_{h},v_{h})+\dfrac{1}{2}a^{n}_{h}(\xi^{n-1}_{h},v_{h})+\dfrac{\gamma_{D}}{h}(\xi^{n}_{h},v_{h})_{\partial \Omega^{n}}\\[8pt]
+\gamma_{g}g^{n}_{h}(\xi^{n}_{h},v_{h})
=\dfrac{1}{2}\mathcal{E}^{n-1}_{f}(v_{h})+\mathcal{E}^{n}_{c}(v_{h})-\mathcal{E}^{n}_{I}(v_{h}),
\end{array}
\end{equation}
where the interpolation error $ \mathcal{E}^{n}_{I}(v_{h}) $ is given  by
\begin{equation}\label{interr}
\mathcal{E}^{n}_{I}(v_{h})=\left(D^{-}_{\Delta t}\eta^{n},v_{h}
\right)_{\Omega^{n}}+\dfrac{1 }{2}a^{n}_{h}(\eta^{n},v_{h})+\dfrac{1}{2}a^{n}_{h}(\eta^{n-1},v_{h})+\dfrac{\gamma_{D}}{h}(\eta^{n},v_{h})_{\partial \Omega^{n}}+\gamma_{g}g^{n}_{h}(\eta^{n},v_{h}).
\end{equation}
We remark that we are using smooth extensions of, e.g., $\eta^{n-1}$, in the following without further notice, whenever variables would be undefined on parts of $\Omega^n$.

\begin{lemma}[Interpolation error]\label{interpollemma}
Let $u \in L^{\infty}(I_{n},H^{m+1}(\Omega(t)), \,u_{t} \in L^{\infty}(I_{n},H^{m}(\Omega(t))$ and $v_{h} \in V^{n,m}_{h}$.
Based on the assumptions~\ref{ass:regularity} and~\ref{ass:CFL}, the interpolation error is bounded by	
\begin{equation}\label{intteq}
\begin{array}{ll}
\Delta t|\mathcal{E}^{n}_{I} (v^{n}_h)| &\leq   c \Delta t \,h^{2m}\mathcal{R}_{I}(u)^2+ \dfrac{\Delta t}{8}\|v^{n}_{h}\|^{2}_{\Omega^{n}}+\dfrac{\Delta t}{32}\|\nabla v^{n}_{h}+\nabla v^{n-1}_{h}\|^{2}_{\Omega^{n}}+\dfrac{1}{32}\| v^{n}_{h}- v^{n-1}_{h}\|^{2}_{\Omega^{n}}\\[8pt]
& \hspace*{3cm}+\dfrac{\Delta t}{16h}\| v^{n}_{h}\|^{2}_{\partial \Omega^{n}}+ \dfrac{\Delta t}{8} g^{n}_{h}(v^{n}_{h},v^{n}_{h}),
\end{array}	
\end{equation}
where $\mathcal{R}_{I}(u)= \|u_{t}\|_{\infty,m}+\|u\|_{\infty,m+1}$.
\end{lemma}
\noindent
{\bf Proof.} The first term in the interpolation error $\mathcal{E}^{n}_{I}(v_{h})$ from \eqref{interr} is estimated as follows:
\begin{equation}\label{Inteq01}
\begin{array}{ll}
\left|\left(\eta^{n}-\eta^{n-1},v^{n}_{h}
\right)_{\Omega^{n}}\right| & \leq  \|\eta^{n}-\eta^{n-1}\|_{\Omega^{n}}\|v^{n}_{h}\|_{\Omega^{n}}\\[10pt]
& = \displaystyle  \left\|\int^{t_{n}}_{t_{n-1}}\eta_{t}(s)ds\right\|_{\Omega^{n}}\|v^{n}_{h}\|_{\Omega^{n}}
\leq \Delta t h^{m}\|\partial_{t}E^{n}u\|_{\infty,m,I_n}\|v^{n}_{h}\|_{\Omega^{n}}.
\end{array}
\end{equation} 
Using the stability of the extension (\ref{exteq2}), we obtain
\begin{equation}\label{Inteq1}
\left|\left(\eta^{n}-\eta^{n-1},v^{n}_{h}
\right)_{\Omega^{n}}\right| \leq c \Delta t \, h^{m} \left(\|u_{t}\|_{\infty,m}+\|u\|_{\infty,m+1}\right)\|v^{n}_{h}\|_{\Omega^{n}}.
\end{equation} 
Next, we use interpolation estimates (\ref{interr1}) and  (\ref{interr2}) in combination with (\ref{exteq1}) to deduce for $k\in \{n-1,n\}$ that
\begin{equation}\label{Inteq20}
\begin{array}{ll}
\dfrac{\Delta t}{2}\left|a^{n}_{h}(\eta^{k},v^{n}_{h})\right|  & \leq  \dfrac{\Delta t}{2} \|\nabla \eta^{k}\|_{\Omega^{n}} \|\nabla v^{n}_{h}\|_{\Omega^{n}}+\Delta t \|\partial_{n} \eta^{k}\|_{\partial \Omega^{n}}\| v^{n}_{h}\|_{\partial \Omega^{n}} \\[12pt]
& \leq c \Delta t\,h^{m} \|u^k\|_{H^{m+1}(\Omega^{k})}\|\nabla v^{n}_{h}\|_{\Omega^{n}}+c \Delta t\,h^{m-1/2} \|u^k\|_{H^{m+1}(\Omega^{k})}\| v^{n}_{h}\|_{\partial \Omega^{n}}.
\end{array}
\end{equation}
Using that, by an inverse inequality and the CFL condition $\Delta t \leq c_{CFL} h^2$,
\begin{align*}
  \|\nabla v_h^n\|_{\Omega^n} &\leq \dfrac{1}{2}\|\nabla v^{n}_{h}+\nabla v^{n-1}_{h}\|_{\Omega^{n}}+\dfrac{1}{2}\|\nabla v^{n}_{h}-\nabla v^{n-1}_{h}\|_{\Omega^{n}} \\ &\leq \dfrac{1}{2}\|\nabla v^{n}_{h}+\nabla v^{n-1}_{h}\|_{\Omega^{n}}+\dfrac{c}{h}\| v^{n}_{h}- v^{n-1}_{h}\|_{\Omega^{n}}\\
  &\leq \dfrac{1}{2}\|\nabla v^{n}_{h}+\nabla v^{n-1}_{h}\|_{\Omega^{n}}+\dfrac{c}{\Delta t^{1/2}}\| v^{n}_{h}- v^{n-1}_{h}\|_{\Omega^{n}}
\end{align*}
we obtain further
\begin{equation}\label{Inteq2}
\begin{array}{ll}
\dfrac{\Delta t}{2}\left|a^{n}_{h}(\eta^{k},v^{n}_{h})\right| 
& \leq  c \Delta t^{1/2} \,h^{m} \|u^k\|_{H^{m+1}(\Omega^{n})}\left(\dfrac{\Delta t^{1/2}}{2}\|\nabla v^{n}_{h}+\nabla v^{n-1}_{h}\|_{\Omega^{n}}+\| v^{n}_{h}- v^{n-1}_{h}\|_{\Omega^{n}}\right)\\[12pt]
& \qquad 
+c \Delta t^{1/2}\,h^{m} \|u^k\|_{H^{m+1}(\Omega^{n})} \dfrac{\Delta t^{1/2}}{h^{1/2}}\| v^{n}_{h}\|_{\partial \Omega^{n}}.
\end{array}
\end{equation}
For the Nitsche penalty term, we have  
\begin{equation}\label{Inteq4}
\begin{array}{ll}
\Delta t \left|\dfrac{\gamma_{D}}{h}(\eta^{n},v^{n}_{h})_{\partial \Omega^{n}}\right| & \leq  c \Delta t \dfrac{\gamma_{D}}{h} \|\eta^{n}\|_{\partial  \Omega^{n}}\| v^{n}_{h}\|_{\partial \Omega^{n}}   \leq c \Delta t^{1/2} \,  h^{m} \|u^{n}\|_{H^{m+1}(\Omega^{n})}  \dfrac{\Delta t^{1/2}}{h^{1/2}} \|v^{n}_{h}\|_{\partial \Omega^{n}}.
\end{array}
\end{equation}
 Finally, the ghost-penalty term is approximated by using $g_h^n(\eta^{n},v^{n}_{h}) \leq g_h^n(\eta^{n},\eta^n)^{1/2} g_h^n(v^{n}_{h}, v_h^n)^{1/2}$ and an interpolation estimate
\begin{equation}\label{Inteq5}
\begin{array}{ll}
\Delta t \left|\gamma_{g}g^{n}_{h}(\eta^{n},v^{n}_{h})\right|  \leq c \Delta t\, h^{m} \|u^{n}\|_{H^{m+1}(\Omega^{n})} g^{n}_{h}(v^{n}_{h},v^{n}_{h})^{1/2}.
\end{array}
\end{equation}
The statement~\eqref{intteq} follows by combining estimates (\ref{Inteq1})-(\ref{Inteq5}) and using Young's inequality. 
\hfill \eop

\subsection{Convergence estimate }

\begin{lemma}[Discrete error]\label{errthm}
	Let $u\in L^{\infty}(I_{n},H^{m+1}(\Omega(t))) \cap W^{1,\infty}(I_{n},H^{m}(\Omega(t))\cap W^{3,\infty}(Q)$
	 be the solution of (\ref{C1}) and $\{u^{k}_{h}\}_{k=1}^n$ the discrete solution of (\ref{2FE1}), respectively. Under  Assumptions~\ref{ass:regularity} and~\ref{ass:CFL}, the discrete error term $\xi^{n}$ satisfies for $\gamma_g, \gamma_D$  sufficiently large
\begin{equation}\label{finalerr}
\|\xi^{n}\|^{2}_{\Omega^{n}}+\Delta t \|\nabla \xi^{n}\|^{2}_{\Omega^{n}} +\Delta t \sum_{k=1}^{n} \mathfrak{E}^k(\xi^{k}_{h}, \xi^{k-1}_{h}) \leq  \exp(c_{T \ref{errthm}} t_{n}) \left(\Delta t^{4}+h^{2m}\right) \mathcal{R}(u)^2,
\end{equation}
where $\mathcal{R}(u)= {\cal R}_C(u) +\mathcal{R}_{I}(u) +\|u(t_n)\|_{W^{3,\infty}(Q)}$,
with ${\cal R}_C$ and ${\cal R}_I$ specified in Lemma~\ref{Tdomainbd} and~\ref{interpollemma}, respectively.
\end{lemma}
\noindent
{\bf Proof.} 
By taking $v_{h}=2\Delta t \xi^{n}_{h}$ in (\ref{Err2}) and using the argumentation from the stability proof (Theorem~\ref{theo.stab}), see~\eqref{2stb5}, we obtain that
\begin{equation}\label{cnveq1}
\begin{array}{ll}
\|\xi^{n}_{h}\|^{2}_{\Omega^{n}}+\dfrac{1}{4}\|\xi^{n}_{h}-\xi^{n-1}_{h}\|^{2}_{\Omega^{n}}+\dfrac{\Delta t}{2}\|\nabla \xi^{n}_{h}\|_{\Omega^{n}}^{2}+\dfrac{\Delta t}{8}\|\nabla \xi^{n}_{h}+\nabla \xi^{n-1}_{h}\|_{\Omega^{n}}^{2}
+\dfrac{3\gamma_{D} \Delta t}{2 h} \|\xi^{n}_{h}\|^{2}_{\partial \Omega^{n}}\\[10pt]
+2\Delta t\gamma_{g}g_{h}(\xi^{n}_{h},\xi^{n}_{h}) \leq \| \xi^{n-1}_{h}\|^{2}_{\Omega^{n-1}}+\dfrac{\Delta t}{2}\|\nabla \xi^{n-1}_{h}\|^{2}_{\Omega^{n-1}}+c \Delta t\| \xi^{n-1}_{h}\|^{2}_{\partial \Omega^{n-1}}\\[10pt]
\hspace{3.5cm}+c \Delta t g^{n-1}_{h}( \xi^{n-1}_{h}, \xi^{n-1}_{h})+ 2\Delta t \left( |\mathcal{E}^{n-1}_{f}(\xi^{n}_{h})|+|\mathcal{E}^{n}_{c}(\xi^{n}_{h})|+|\mathcal{E}^{n}_{I}(\xi^{n}_{h})|\right).
\end{array}
\end{equation}
By combining results from Lemmas~\ref{Tdomainbd},~\ref{cnstlemma} and~\ref{interpollemma} we have
\begin{align}\label{rhss}
\begin{split}
 2\Delta t &\left( |\mathcal{E}^{n-1}_{f}(\xi^{n}_{h})|+|\mathcal{E}^{n}_{c}(\xi^{n}_{h})|+|\mathcal{E}^{n}_{I}(\xi^{n}_{h})|\right) \\[12pt]
 &\qquad\leq c \Delta  t \left(\Delta t^{4}+h^{2m}\right)  \mathcal{R}(u)
+ \dfrac{\Delta t}{4}\|\xi^{n}_{h}\|^{2}_{\Omega^{n}}+\dfrac{\Delta t}{16}\|\nabla \xi^{n}_{h}+\nabla \xi^{n-1}_{h}\|^{2}_{\Omega^{n}}+\dfrac{1}{16}\| \xi^{n}_{h}- \xi^{n-1}_{h}\|^{2}_{\Omega^{n}}\\[12pt]
& \qquad\qquad + \dfrac{\Delta t}{8h} \|\xi^{n}_{h}\|^{2}_{\partial \Omega^{n}}+c \Delta t g^{n}_{h}(\xi^{n}_{h},\xi^{n}_{h})+c \Delta t g^{n-1}_{h}(\xi^{n-1}_{h},\xi^{n-1}_{h}).
\end{split}
\end{align}
Inserting~\eqref{rhss} into~\eqref{cnveq1} and absorbing terms into the left-hand side yields for $n>1$ and $\gamma_D, \gamma_g$ sufficiently large 
\begin{align}\label{ferreq3}
\begin{split}
\left(1-\dfrac{\Delta t}{4}\right)\|&\xi^{n}_{h}\|^{2}_{\Omega^{n}}+\dfrac{1}{4}\|\xi^{n}_{h}-\xi^{n-1}_{h}\|^{2}_{\Omega^{n}}+\dfrac{\Delta t}{4}\|\nabla \xi^{n}_{h}\|_{\Omega^{n}}^{2}+\dfrac{\Delta t}{16}\|\nabla \xi^{n}_{h}+\nabla \xi^{n-1}_{h}\|_{\Omega^{n}}^{2}\\[10pt]
&+\dfrac{\gamma_{D} \Delta t}{h} \|\xi^{n}_{h}\|^{2}_{\partial \Omega^{n}} +\Delta t\gamma_{g}g_{h}(\xi^{n}_{h},\xi^{n}_{h})
\leq \| \xi^{n-1}_{h}\|^{2}_{\Omega^{n-1}}+\dfrac{\Delta t}{2}\|\nabla \xi^{n-1}_{h}\|^{2}_{\Omega^{n-1}}+c \Delta t\| \xi^{n-1}_{h}\|^{2}_{\partial \Omega^{n-1}}\\[10pt]
&\hspace{6cm}+c \Delta t g^{n-1}_{h}( \xi^{n-1}_{h}, \xi^{n-1}_{h})+ c \Delta  t \left(\Delta t^{4}+h^{2m}\right)  \mathcal{R}(u).
\end{split}
\end{align}

As $\xi_h^0=0$, we obtain for $n=1$
\begin{align*}
\begin{split}
\left(1-\dfrac{\Delta t}{4}\right)\|&\xi^{1}_{h}\|^{2}_{\Omega^{1}}+\dfrac{1}{4}\|\xi^{1}_{h}-\xi^{0}_{h}\|^{2}_{\Omega^{1}}+\dfrac{\Delta t}{4}\|\nabla \xi^{1}_{h}\|_{\Omega^{1}}^{2}+\dfrac{\Delta t}{16}\|\nabla \xi^{1}_{h}+\nabla \xi^{0}_{h}\|_{\Omega^{1}}^{2}\\[10pt]
&+\dfrac{\gamma_{D} \Delta t}{h} \|\xi^{1}_{h}\|^{2}_{\partial \Omega^{1}} +\Delta t\gamma_{g}g_{h}^1(\xi^{1}_{h},\xi^{1}_{h}) \,\leq\, c \Delta  t \left(\Delta t^{4}+h^{2m}\right)  \mathcal{R}(u).
\end{split}
\end{align*}
Now, summing over $k=1,2,\ldots,n$, we deduce that
\begin{equation}\label{ferreq4}
\begin{array}{ll}
\displaystyle
\|\xi^{n}_{h}\|^{2}_{\Omega^{n}}+\Delta t\|\nabla \xi^{n}_{h}\|_{\Omega^{n}}^{2}+\dfrac{\Delta t}{8}\sum_{k=1}^{n}\mathfrak{E}^k(\xi^{k}_{h}, \xi^{k-1}_{h}) \, \leq\,  c \, t_{n} \left(\Delta t^{4}+h^{2m}\right) \mathcal{R}(u)+c \Delta t \sum_{k=2}^{n}\| \xi^{k-1}_{h}\|^{2}_{\Omega^{k-1}}.
\end{array}
\end{equation}
Finally, we use the  discrete Gronwall lemma  to conclude the result. This completes the proof.
\hfill 
\eop
\begin{theorem}[Global error]\label{globalerrthm} 
Based on the assumption made in Lemma~\ref{errthm}, the global error $e^k=u^k-u_h^k, k=1,...,n$ satisfies
\begin{equation}\label{globalerr}
\|e^{n}\|^{2}_{\Omega^{n}}+\Delta t \|\nabla e^{n}\|^{2}_{\Omega^{n}} +\Delta t \sum_{k=1}^{n} \mathfrak{E}^k(e^{k}, e^{k-1})  \leq  \exp(c_{T \ref{globalerrthm}} t_{n}) \left(\Delta t^{4}+h^{2m}\right) \mathcal{R}(u)^2,
\end{equation}
where ${\cal R}(u)$ is defined in Lemma~\ref{errthm}.
\end{theorem}
\noindent
{\bf Proof.} Using the interpolation error estimates (\ref{interr1}), (\ref{interr2}) and (\ref{Inteq1}), we deduce that 
\begin{equation}\label{finterreq}
\|\eta^{n}_{h}\|^{2}_{\Omega^{n}}+\Delta t\|\nabla \eta^{n}_{h}\|_{\Omega^{n}}^{2}+\Delta t \sum_{k=1}^{n} \mathfrak{E}^k(\eta^{k}_{h}, \eta^{k-1}_{h}) \leq  c \, h^{2m}\mathcal{R}_{I}(u)^2.
\end{equation}
In combination with Lemma~\ref{errthm}, this proves~\eqref{globalerr}.
\hfill \eop

\section{Numerical results}\label{Sec:5}
In this section, we show numerical results in two and three space dimensions to  verify the theoretical findings and the practical behaviour of the numerical method. All the numerical experiments have been obtained using the CutFEM
library~\cite{burman_ijnme15}, which is based on FEniCS \cite{alnaes2015fenics}.

To verify the theoretical results, we will analyse the error terms $e^{k}=u^k-u^{k}_{h}, k=1,...,n$ in the following norms
\begin{align*}
\|e^n\|_{L^2(\Omega^n)},\,\, \|e\|_{L^{2}(L^{2})} &=\left(\Delta t \sum_{k=1}^{n}\|e^{k}\|^{2}_{L^{2}(\Omega^{k})}\right)^{1/2},\,\, 
 \|e\|_{L^{2}(H_{av}^{1})}&=\left(\Delta t \sum_{k=1}^{n}\|\nabla e^{k} + \nabla e^{k-1}\|^{2}_{\Omega^{k}}\right)^{1/2}.
\end{align*}
Given the respective CFL condition (Assumption~\ref{ass:CFL}), Theorem~\ref{globalerrthm} guarantees second-order convergence in time and convergence of order $m$ in space in the $L^2$-norm at the end time $\|e^n\|_{L^2(\Omega^n)}$ and in the averaged $L^2(H^1)$-norm $\|e\|_{L^{2}(H_{av}^{1})}$.


\subsection{2d example}
\label{sec:51}

\begin{example}\label{Exmp1}
We consider a circle traveling with constant velocity $w=(1,0)$ towards the right. The domain is given by
\[
\Omega(t)=\{(x,y): \,(x-0.5-t)^{2}+(y-0.5)^{2} \leq 0.9\}
\]
in the time interval $I = [0, 0.1]$. The data of the model example is chosen  in such a way that the exact solution is 
\[
u(x,y,t)=\exp(-4\pi^{2}t)\cos(2\pi x) \cos(2\pi y).
\]
\end{example}

\begin{figure}[t]
\centerline{
\begin{tabular}{ccc}
\resizebox*{5cm}{!}{\includegraphics{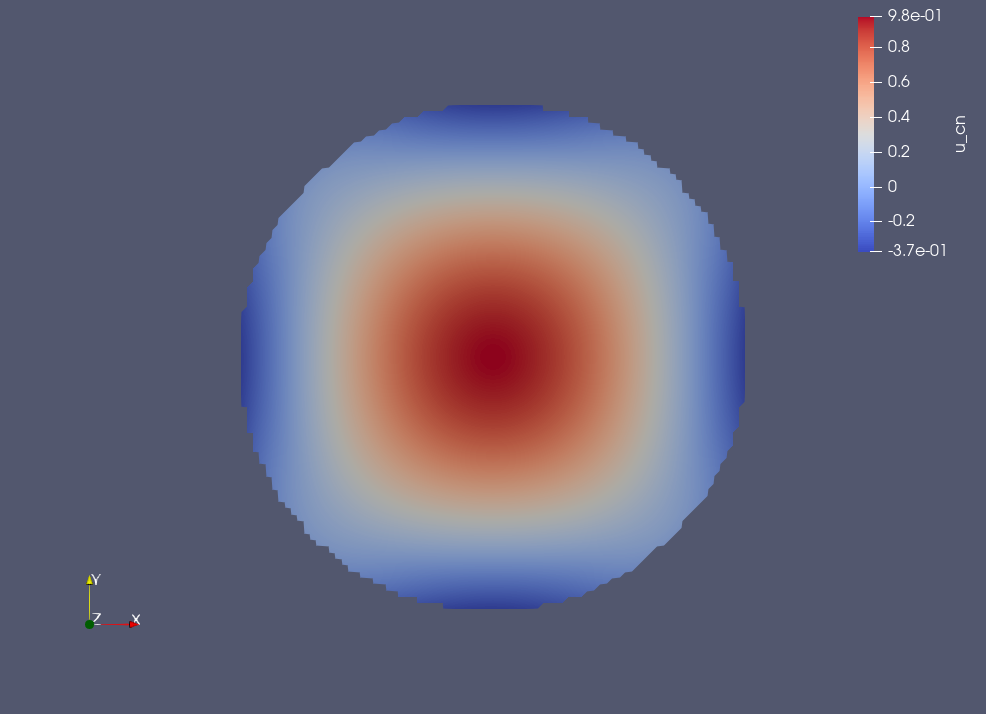}}%
& 
\resizebox*{5cm}{!}{\includegraphics{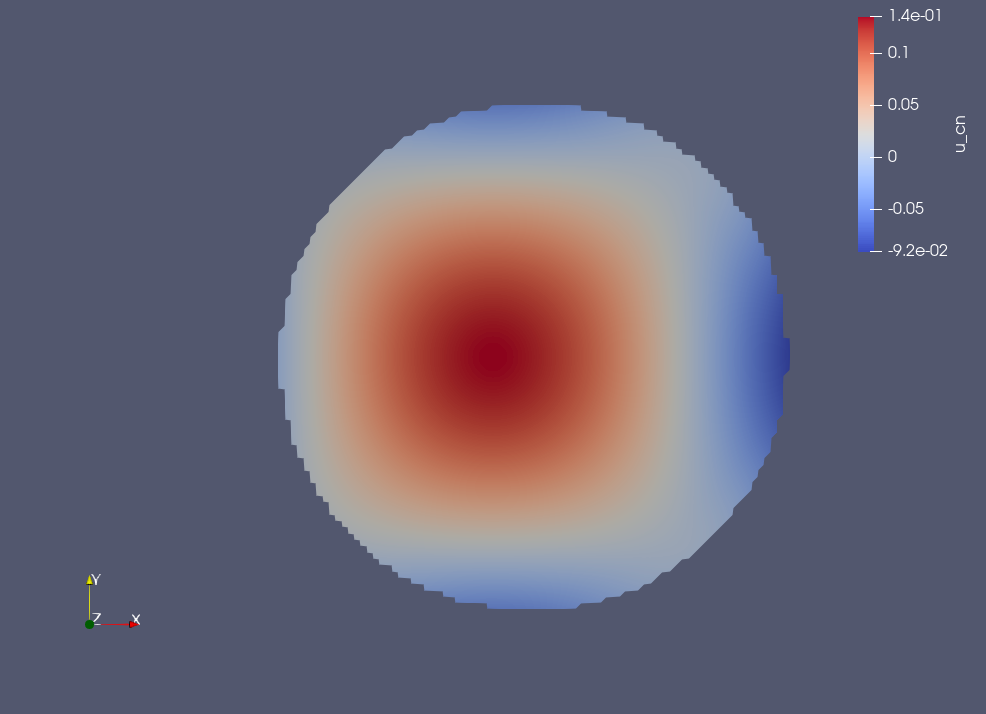}}%
& 
\resizebox*{5cm}{!}{\includegraphics{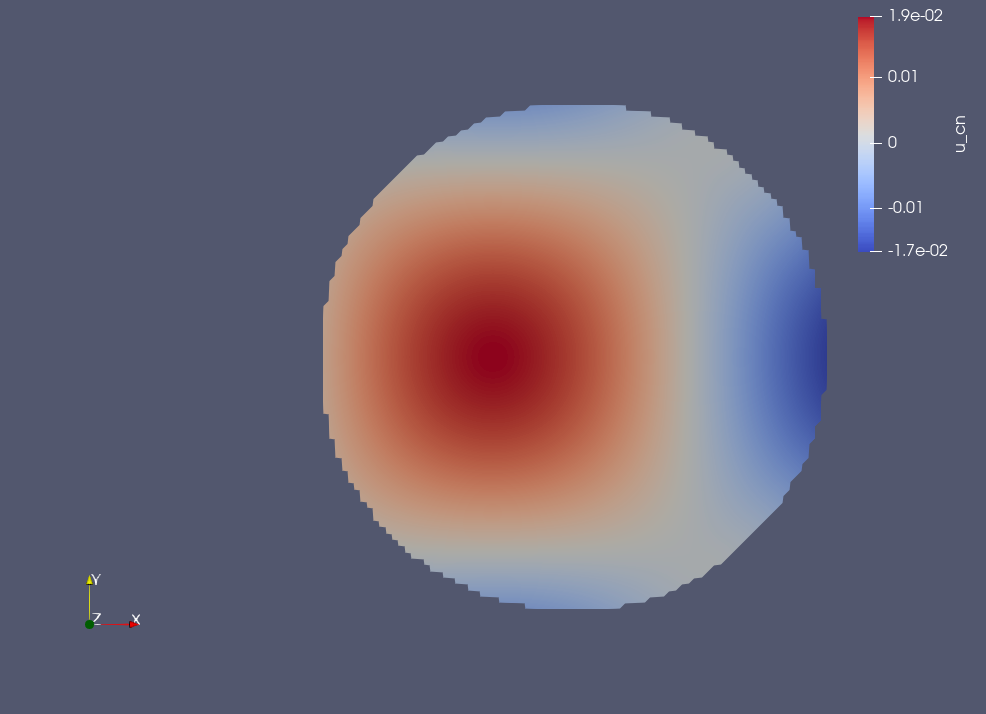}}%
\end{tabular}
} 
\caption{\label{fig:numres2d} 
Illustration of the numerical solution $u$ at $t=0$ (left), $t=0.05$ (center) and $t=0.1$ (right).}
\end{figure}

An illustration of the numerical solution is given in Figure~\ref{fig:numres2d}. As background domain $D$, we use the unit square, i.e.\,$D=[0,1]^{2}$. The  background 
triangulations $\mathcal{T}_h$ are created by a uniform subdivision of the unit square into triangles and successive refinement. For each time step $n$, the active triangulation $\mathcal{T}^{n}_{h,\delta}$ is then extracted from ${\cal T}_h$, as described in Section~\ref{Sec:2}. We use $\delta =4\Delta t$.
 
 \paragraph{Estimated orders of convergence}
We will show results for different time-step sizes $\Delta t_i = \nicefrac{1}{50} \cdot 2^{-i}, i=0,...,4$ and mesh sizes $h_j=\nicefrac{1}{32} \cdot 2^{-j}, j=0,...,3$. From the computed errors, we will estimate the temporal and spatial order of convergence. Therefore, we assume that the total error can be decomposed  into a temporal and a spatial component as follows
\[
g(\Delta t, h) = g_{\Delta t}(\Delta t)+g_{h}(h) = c_{\Delta t}\Delta t^{eoc_{\Delta t}} + c_{h}h^{eoc_{h}}
\]
with constants $c_{\Delta t}, c_h$ and estimated orders of convergence $eoc_{\Delta t},\, eoc_h$. For a fixed mesh size $h_j$, this relation becomes
\begin{align}\label{temporal}
g(\Delta t, \cdot)  = g_{h_j}+ c_{\Delta t} \Delta t^{eoc_{\Delta t}}
\end{align}
with a fixed spatial error part $g_{h_j}$. We will use~\eqref{temporal} to estimate the order of convergence in time $eoc_{\Delta t}$ by means of a least-squares fit of all available error values for a fixed $h_j$ to find the three parameters $g_{h_j}, \,c_{\Delta t}$ and $eoc_{\Delta t}$. 
Analogously, we estimate the spatial order of convergence $eoc_h$ by a least-squares fit of the function
\begin{align}\label{spatial}
g(\cdot,h)  = g_{\Delta t_i}+ c_h h^{eoc_{h}}
\end{align}
against all available error values for a fixed time-step size $\Delta t_i$.

Finally, we will also compute estimated orders of convergence for the ''diagonal values'', which correspond to fixing $\Delta t =\bar{c}h$ for $\bar{c}\in \{\nicefrac{32}{50}, \nicefrac{32}{100}\}$. Here we fit the two parameters $c_{\Delta t, h}$ and $eoc_{\Delta t,h}$ of the function
\begin{align}\label{spacetime}
g(\bar{c}h,h)  =  c_{\Delta t, h} h^{eoc_{\Delta t,h}}.
\end{align}
against the computed error values.



\subsubsection{$P_1$ finite elements}
Firstly, we consider $P_1$ finite elements. We choose the numerical parameters as $\gamma_D=1$ and $\gamma_g=10^{-3}$. The errors in the $L^{2}$-norm at the end time, the $L^{2}(L^{2})$-norm and the $L^{2}(H_{av}^{1})$-norm are shown in Table~\ref{tabP1} for different $\Delta t$ and $h$. The estimated orders of convergence are shown, if the asymptotic standard error (computed by gnuplot~\cite{williams2017gnuplot}) was below 20\%; otherwise we draw a '-'.

\begin{table}[t]
		\centering
		End-time error $\|e^n\|_{L^2(\Omega^n)}$
		\label{table2dP1L2}
		\begin{tabular}{c |c c c c c | c } 
			\hline
			$h   \downarrow /\Delta t \rightarrow$  & $\nicefrac{1}{50}$ & $\nicefrac{1}{100}$ & $\nicefrac{1}{200}$ & $\nicefrac{1}{400}$ & $\nicefrac{1}{800}$    & \textbf{eoc}$_{\Delta t}$\\ [0.5ex]
			\hline
	    	    $\nicefrac{1}{32}$ &\dashuline{1.93e-03}	& \underline{5.40e-04}	& 3.14e-04 &	 2.78e-04	& 2.68e-04	&\textbf{2.62}\\
                $\nicefrac{1}{64}$  &	1.21e-03	& \dashuline{2.77e-04}	& \underline{1.08e-04} &	 7.84e-05 &	 7.17e-05 &	\textbf{2.46}\\
                $\nicefrac{1}{128}$ & 8.88e-04  &	 1.96e-04 &	 \dashuline{4.95e-05} &	 \underline{2.54e-05}	& 1.96e-05	&\textbf{2.29}\\
                $\nicefrac{1}{256}$ &	7.78e-04 &	 1.57e-04 &	 3.18e-05 &	 \dashuline{1.15e-05} &	 \underline{6.18e-06} &\textbf{2.35}\\
			\hline
			\textbf{eoc}$_{h}$  & 	\textbf{1.26} &		\textbf{1.55}	&	\textbf{1.80}	&	\textbf{1.91}		& \textbf{1.92}& \\
			\textbf{eoc}$_{\Delta t, h}$  & 	 &	&	&	\dashuline{\textbf{2.79}}		& \underline{\textbf{2.30}}&\\
			\hline
		\end{tabular}
\bigskip

		Error $\|e\|_{L^{2}(L^{2})}$\\
		\label{table2dP1L2L2}
		\begin{tabular}{c |c c c c c | c } 
			\hline
			$h   \downarrow /\Delta t \rightarrow$  & $\nicefrac{1}{50}$ & $\nicefrac{1}{100}$ & $\nicefrac{1}{200}$ & $\nicefrac{1}{400}$ & $\nicefrac{1}{800}$    & \textbf{eoc}$_{\Delta t}$\\ [0.5ex]
			\hline
	    	    $\nicefrac{1}{32}$ & \dashuline{1.44e-03}	& \underline{8.42e-04} & 7.12e-04 & 	 6.87e-04	& 6.83e-04	&\textbf{2.25}\\
                $\nicefrac{1}{64}$  &	9.69e-04 	& \dashuline{3.62e-04}	& \underline{2.21e-04} &	1.89e-04 & 	 1.82e-04 &	\textbf{2.11}\\
                $\nicefrac{1}{128}$ & 8.35e-04  &	2.37e-04	& \dashuline{9.13e-05} &	 \underline{5.72e-05} &	 4.89e-05	&\textbf{2.04}\\
                $\nicefrac{1}{256}$ &	7.97e-04 & 2.01e-04 &	 5.75e-05 &	 \dashuline{2.30e-05}	& \underline{1.45e-05} &\textbf{2.05}\\
			\hline
			\textbf{eoc}$_{h}$  & 	\textbf{1.81} &		\textbf{1.91}	&	\textbf{1.92}	&	\textbf{1.92}		& \textbf{1.92}& \\
				\textbf{eoc}$_{\Delta t, h}$  & 	 &	&	&	\dashuline{\textbf{1.99}}		& \underline{\textbf{1.93}}&\\
			\hline
		\end{tabular}
		\bigskip

		Error $\|e\|_{L^{2}(H_{av}^{1})}$\\
		\label{table2dP1L2H1}
		\begin{tabular}{c |c c c c c | c } 
			\hline
			$h   \downarrow /\Delta t \rightarrow$  & $\nicefrac{1}{50}$ & $\nicefrac{1}{100}$ & $\nicefrac{1}{200}$ & $\nicefrac{1}{400}$ & $\nicefrac{1}{800}$    & \textbf{eoc}$_{\Delta t}$\\ [0.5ex]
			\hline
	    	    $\nicefrac{1}{32}$ & \dashuline{4.77e-02} &	 \underline{4.53e-02} &	 4.50e-02	& 4.49e-02 &	 4.49e-02	&\textbf{2.83}\\
                $\nicefrac{1}{64}$  &	2.64e-02 &	\dashuline{2.36e-02} &	 \underline{2.33e-02} &	2.33e-02 & 	 2.33e-02 &	\textbf{3.33}\\
                $\nicefrac{1}{128}$ & 1.65e-02 &	 1.23e-02 &	\dashuline{1.19e-02} & 	 \underline{1.19e-02} &	 1.19e-02	&\textbf{3.32}\\
                $\nicefrac{1}{256}$ &	1.26e-02 &	6.72e-03 & 	 6.06e-03 &	 \dashuline{6.01e-03}	& \underline{6.00e-03}  &\textbf{3.21}\\
			\hline
			\textbf{eoc}$_{h}$  & 	\textbf{1.17} &		\textbf{0.96}	&	\textbf{0.93}	&	\textbf{0.93}		& \textbf{0.93}& \\
			\textbf{eoc}$_{\Delta t, h}$  & 	 &	&	&	\dashuline{\textbf{1.00}}		& \underline{\textbf{0.96}}&\\
			\hline
		\end{tabular}
		\caption{$L^2(T), L^2(L^2)$ and $L^2(H^1_{\text{av}})$-norm errors for $P_1$ finite elements applied to Example 5.1. The estimated orders of convergence are computed according to~\eqref{temporal}-\eqref{spacetime}. The diagonal orders are computed from the underlined error values.\label{tabP1}}
\end{table}

		We observe estimated spatial convergence orders close to two in the $L^2$-norms and close to one in the $L^2(H_{\text{av}}^1)$-norm. Note that Theorem~\ref{globalerrthm} guarantees only first-order convergence in space. We expect, however, that using a duality argument second-order convergence in space could be shown in the $L^2(L^2)$-norm, as in~\cite{burman2022eulerian}. 
		
		The estimated temporal orders of convergence are close to two or larger in the $L^2$-norms. This is in agreement with Theorem~\ref{globalerrthm}. In the $L^2(H_{\text{av}}^1)$-norm the estimated eoc$_{\Delta t}$ seems to be even larger than three. This has to be read carefully, however, as the spatial error part clearly dominates the overall error in this case.
		
		Finally, the diagonal orders are around two in the $L^2(L^2)$-norm and even slightly higher in the $L^2$-norm at the end time. In the $L^2(H_{\text{av}})$-norm the spatial part is dominant and we obtain eoc$_{\Delta t, h}$ close the one, in agreement with Theorem~\ref{globalerrthm}.

\subsubsection{$P_2$ finite elements}
Next, we consider Example~\ref{Exmp1} with $P_2$ finite elements.
We increase the Nitsche parameter to $\gamma_D=10$, as for higher polynomial degree a larger Nitsche parameter is required, see e.g.\,~\cite{Johanssonetal2015}. The ghost-penalty parameter $\gamma_g$ is still chosen as $10^{-3}$, but now, according to (\ref{ghostdef}), second derivatives are included in the ghost-penalty term. The $L^{2}$-norm errors at the end time, the $L^{2}(L^{2})$- and the $L^2(H_{av}^1)$-norm errors are shown in Table~\ref{tabP2}.

\begin{table}[t]
		\center
		End-time error $\|e^n\|_{L^2(\Omega^n)}$\\
		\label{table2dP2L2}
		\begin{tabular}{c |c c c c c | c } 
			\hline
			$h   \downarrow /\Delta t \rightarrow$  &  $\nicefrac{1}{100}$ & $\nicefrac{1}{200}$ & $\nicefrac{1}{400}$ & $\nicefrac{1}{800}$    & $\nicefrac{1}{1\,600}$ & \textbf{eoc}$_{\Delta t}$\\ [0.5ex]
			\hline
	    	    $\nicefrac{1}{32}$ &  \dashuline{1.16e-04} & \underline{3.09e-05} & 1.08e-05	& 5.80e-06 &	4.66e-06	&\textbf{2.06}\\
                $\nicefrac{1}{64}$  &	1.17e-04	& \dashuline{2.66e-05}	& \underline{7.33e-06} &         2.78e-06 &	1.44e-06 &	\textbf{2.18}\\
                $\nicefrac{1}{128}$ & 1.11e-04 	&  2.43e-05     &   \dashuline{6.23e-06} &        \underline{1.76e-06} & 	6.78e-07	&\textbf{2.22}\\
                $\nicefrac{1}{256}$ &	1.12e-04	& 2.31e-05	& 5.78e-06 &        \dashuline{1.54e-06} &	\underline{4.36e-07} &\textbf{2.29}\\
			\hline
			\textbf{eoc}$_{h}$  & 	\textbf{--} &		\textbf{0.91}	&	\textbf{1.60}	&	\textbf{1.68}		& \textbf{2.01}& \\
				\textbf{eoc}$_{\Delta t, h}$  & 	 &	&	&\dashuline{\textbf{2.12}}		& \underline{\textbf{2.08}}& \\
			\hline
		\end{tabular}
	\bigskip

		Error $\|e\|_{L^{2}(L^{2})}$\\
		\label{table2dP2L2L2}
		\begin{tabular}{c |c c c c c | c } 
			\hline
			$h   \downarrow /\Delta t \rightarrow$  &  $\nicefrac{1}{100}$ & $\nicefrac{1}{200}$ & $\nicefrac{1}{400}$ & $\nicefrac{1}{800}$    & $\nicefrac{1}{1\,600}$ & \textbf{eoc}$_{\Delta t}$\\ [0.5ex]
			\hline
	    	    $\nicefrac{1}{32}$ & \dashuline{1.89e-04} &	 \underline{5.35e-05} &	 2.13e-05 &  	1.27e-05 & 	1.10e-05	&\textbf{2.04}\\
                $\nicefrac{1}{64}$  &	1.84e-04	& \dashuline{4.66e-05}	& \underline{1.36e-05} &	5.09e-06 & 	3.16e-06  &	\textbf{2.03}\\
                $\nicefrac{1}{128}$ & 1.82e-04	& 4.50e-05	& \dashuline{1.16e-05} &   	\underline{3.32e-06} &        1.28e-06	&\textbf{2.03}\\
                $\nicefrac{1}{256}$ &	1.82e-04 &	 4.45e-05 &	 1.12e-05 &	\dashuline{2.88e-06} &	\underline{8.32e-07} &\textbf{2.03}\\
			\hline
			\textbf{eoc}$_{h}$  & 	\textbf{1.72} &		\textbf{2.00}	&	\textbf{1.99}	&	\textbf{2.09}		& \textbf{2.05}& \\
				\textbf{eoc}$_{\Delta t,h}$  & 	&	&	&	\dashuline{\textbf{2.02}}		& \underline{\textbf{1.98}}& \\
			\hline
		\end{tabular}
\bigskip 

		Error $\|e\|_{L^{2}(H_{av}^{1})}$ \\
		\label{table2dP2L2H1}
		\begin{tabular}{c |c c c c c | c } 
			\hline
			$h   \downarrow /\Delta t \rightarrow$  &  $\nicefrac{1}{100}$ & $\nicefrac{1}{200}$ & $\nicefrac{1}{400}$ & $\nicefrac{1}{800}$    & $\nicefrac{1}{1\,600}$ & \textbf{eoc}$_{\Delta t}$\\ [0.5ex]
			\hline
	    	    $\nicefrac{1}{32}$ & \dashuline{3.42e-03}	& \underline{2.25e-03} &	2.17e-03  &	 2.17e-03    &    2.17e-03	&\textbf{3.88}\\
                $\nicefrac{1}{64}$  &	2.81e-03 &	\dashuline{9.11e-04}	& \underline{6.22e-04} &         6.03e-04     &   6.03e-04 &	\textbf{2.82}\\
                $\nicefrac{1}{128}$ & 2.80e-03 &	7.26e-04	& \dashuline{2.45e-04} &          \underline{1.77e-04}  &	1.72e-04  	&\textbf{2.22}\\
                $\nicefrac{1}{256}$ &	2.81e-03 	& 7.14e-04 &	1.87e-04 &           \dashuline{6.84e-05} &	\underline{5.30e-05}   &\textbf{2.03}\\
			\hline
			\textbf{eoc}$_{h}$  & 	\textbf{--} &		\textbf{2.91}	&	\textbf{2.11}	&	\textbf{1.89}		& \textbf{1.86}& \\
				\textbf{eoc}$_{\Delta t, h}$  & 	 &	&	&	\dashuline{\textbf{1.90}}		& \underline{\textbf{1.85}}&\\
			\hline
		\end{tabular}
			\caption{$L^2(T), L^2(L^2)$ and $L^2(H^1_{\text{av}})$-norm errors for $P_2$ finite elements applied to Example 5.1.\label{tabP2}}
\end{table}

	Firstly, we observe that the absolute values of the errors are significantly smaller compared to $P_1$ finite elements. 
	The spatial orders of convergence are close to two in all norms for smaller $\Delta t$. 
	We note that in Theorem~\ref{globalerrthm} second order in space has been shown for the end-time $L^2$-norm and the $L^2(H_{av}^1)$-norm. Using a duality argument, one could even hope for convergence order three in the $L^2(L^2)$-norm. We need to consider, however, that for these results the quadrature error related to a curved boundary has not been taken into account. In the CutFEM library used here, the geometry is approximated linearly in the set-up of the quadrature rule, see~\cite{burman_ijnme15}. This can lead to a reduced order of convergence, namely order 1.5 in the $H^1$-norms and 2 in the $L^2$-norms, see \cite{burman2018cut} for results for a CutFEM approach applied to an elliptic problem on curved domains. This reduction can also be observed in the $L^2$-norm errors in the example considered here.
	
	The estimated temporal orders of convergence are again close to two in the $L^2$-norms, which confirm the estimates in Theorem~\ref{globalerrthm}. In the $L^2(H_{\text{av}}^1)$-norm the error is still clearly dominated by the spatial part for $h\geq \frac{1}{64}$. This changes, however, on the finer levels, where the temporal part gets dominant and the eoc$_{\Delta t}$ is very close to two, in agreement with Theorem~\ref{globalerrthm}. 
	
	The spatial and temporal convergence orders are confirmed by the ''diagonal'' orders, which are close to two in all cases.

	\subsection{3d example}
	
	\begin{example}\label{Exmp2}
We consider  a 3-dimensional rectangular channel with a moving upper and lower wall in the time interval  $I = [0,1]$, inspired by a pump. The moving domain is given by
\[
\Omega(t)=(0,4) \times (-1+0.1\sin t , 1-0.1\sin t) \times (-1,1).
\]
The source term and boundary data is chosen  in such a way that the exact solution of the model problem  (\ref{C1}) is 
\[
u(x,y,z,t)=\exp(-t)\left((1-0.1\sin t)^{2}-y^{2}\right).
\]
\end{example}

\begin{figure}[t]
\centerline{
\begin{tabular}{cc}
\resizebox*{7cm}{!}{\includegraphics{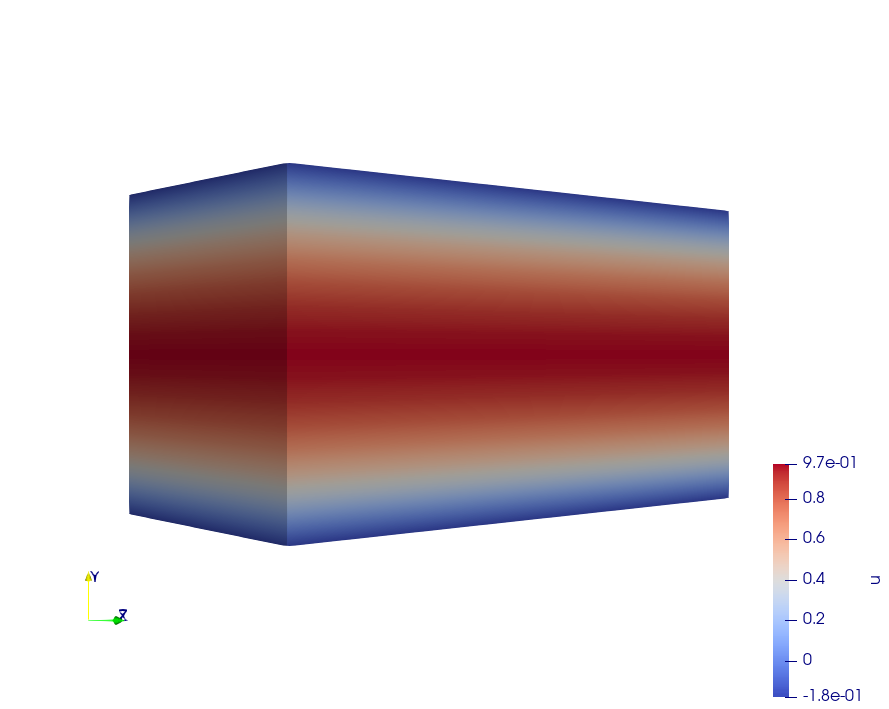}}%
& 
\resizebox*{7cm}{!}{\includegraphics{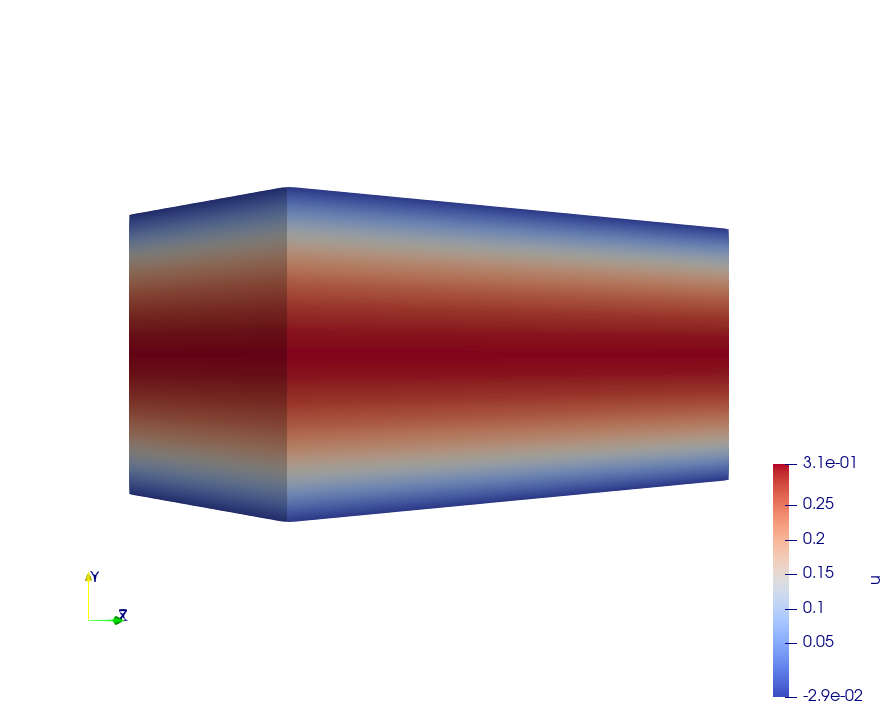}}%
\end{tabular}
} 
\caption{\label{fig:3d} Illustration of the numerical solution of Example 5.2 at time  $t=0$ (left) and $t=1$ (right).}
\end{figure}

As background domain $D$, we use the box $[0,4] \times[-1.1,1.1] \times [-1,1]$. The  background 
triangulations $\mathcal{T}_h$ are created by uniform subdivisions of $D$  into tetrahedra and successive refinement. We use again $\delta =4\Delta t$ and choose $\gamma_{g}=0.1$ for $P_1$ and $\gamma_{g}=1$ for $P_2$ finite elements, respectively and $\gamma_D=10$ in both cases. We note that in this example the quadrature error is zero, as the boundary $\partial\Omega^k$ consists of plane surfaces for all $k$. An illustration of the numerical solution at times $t=0$ and $t=1$  is given in Figure~\ref{fig:3d}.

\subsubsection{$P_1$ finite elements}

\begin{figure}[t]
\centerline{
\begin{tabular}{cc}
\resizebox*{13cm}{!}{\includegraphics{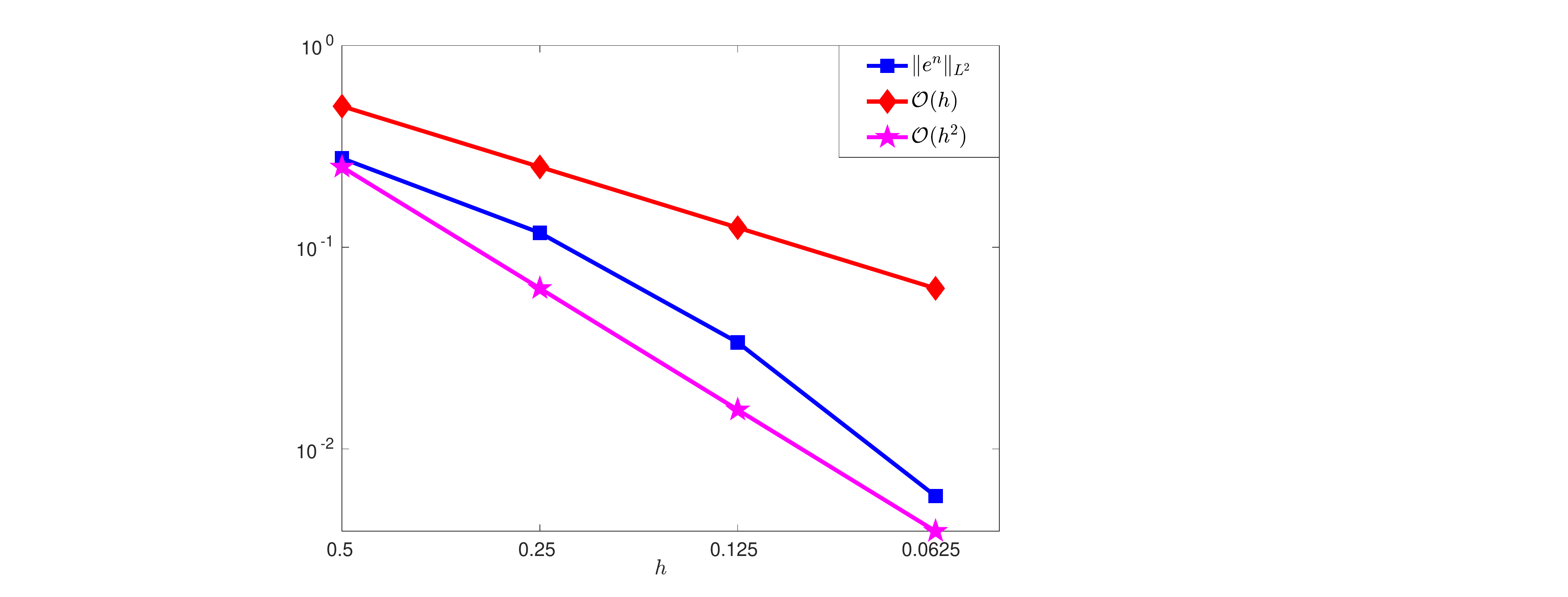}}%
& \hspace{-9cm}
\resizebox*{12.5cm}{!}{\includegraphics{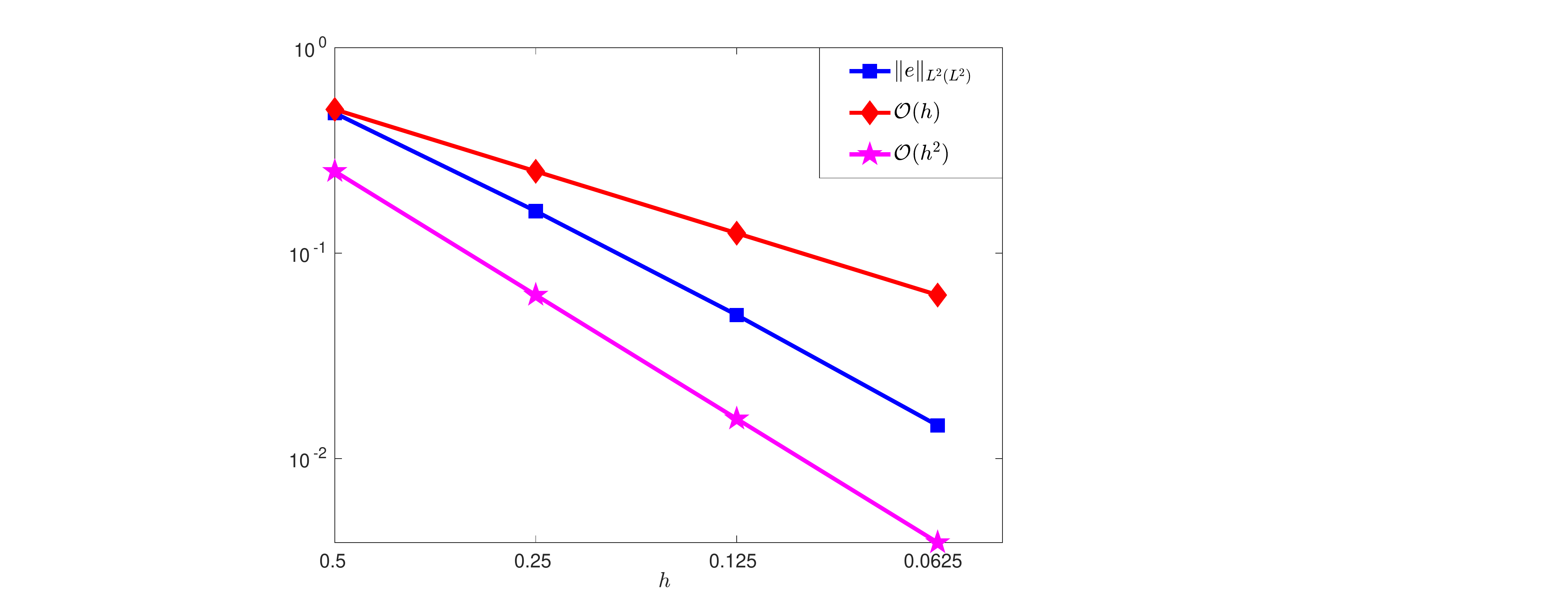}}\\
\hspace*{7.5cm}
\resizebox*{13cm}{!}{\includegraphics{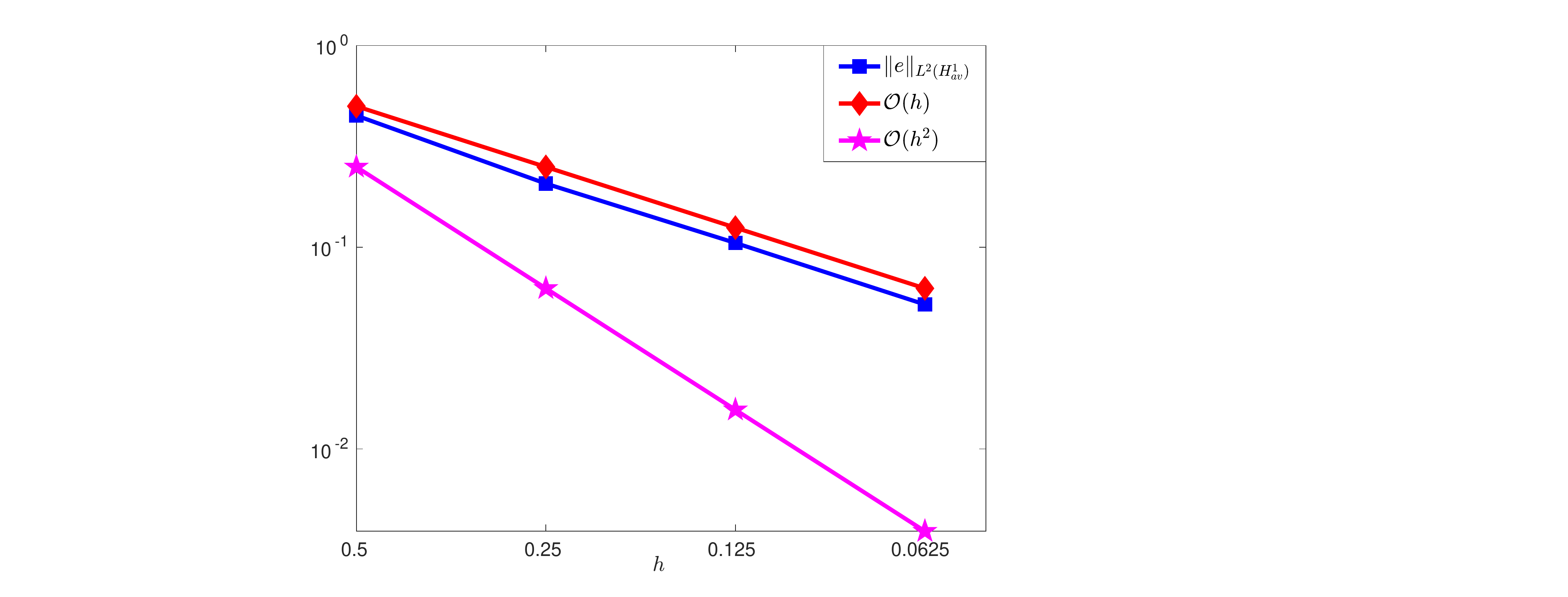}}%
\end{tabular}
} 
\caption{\label{fig:P13d} Errors for example 5.2 (3d)  for $P_1$ finite elements under simultaneous refinement in space and time ($\Delta t = h/10$) for $P_1$ finite elements.  \textbf{Top left}: $L^2$-norm at the end time. \textbf{Top right}: $L^2(L^2)$-norm. \textbf{Bottom}: $L^2(H_{\text{av}}^1)$-norm.}
\end{figure}

As the numerical experiments in three space dimensions are much more time-consuming compared to two dimensions, we focus on simultaneous refinement in space and time by choosing $h_i = \frac{\Delta t_i}{10}=2^{-i-1}$ for $i=0,...,3$. The resulting errors in the four norms introduced above are plotted in Figure~\ref{fig:P13d} over the mesh size (blue curves) and compared to linear (red) and quadratic convergence (pink). We observe second-order convergence in the $L^2$-norms and first-order convergence in the $L^2(H_{\text{av}}^1)$-norm. We note that for $m=1$ first-order convergence in space has been shown in Theorem~\ref{globalerrthm} in the $L^2$-norm at the end time and in the $L^2(H^1_{av})$-norm. The numerical results indicate again that second-order convergence in space could be shown in the $L^2(L^2)$-norm using a duality argument. The first-order convergence in the $L^2(H_{\text{av}}^1)$-norm is optimal, as the spatial error part dominates the overall error.

\subsubsection{$P_2$ finite elements}

\begin{figure}[t]
\centerline{
\begin{tabular}{cc}
\hspace{1cm}
\resizebox*{13cm}{!}{\includegraphics{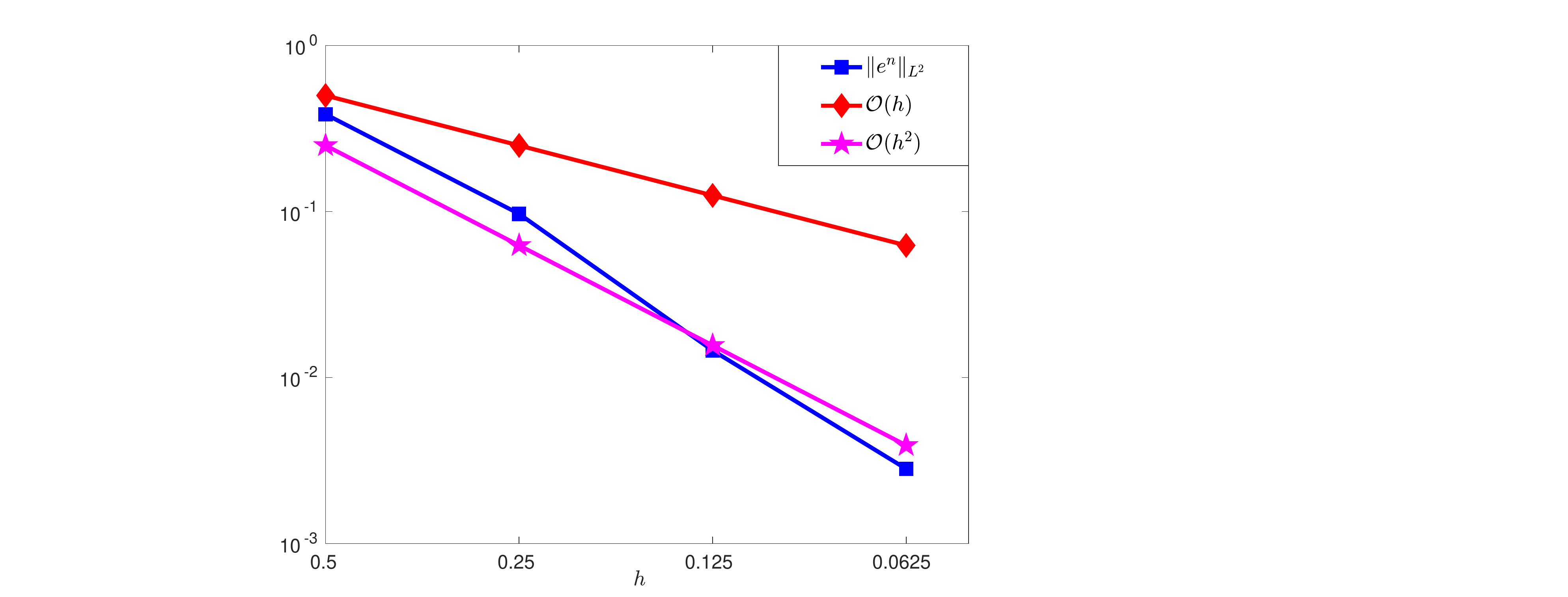}}%
& \hspace{-9cm}
\resizebox*{13.5cm}{!}{\includegraphics{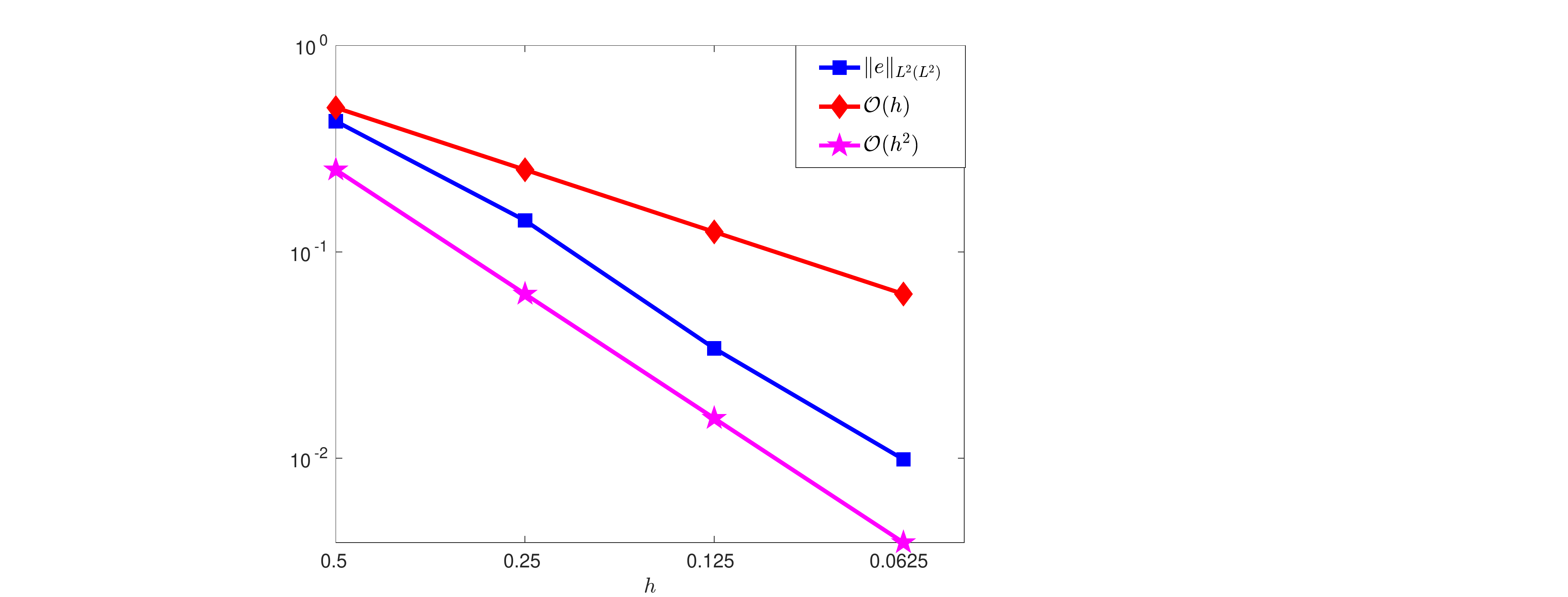}}\\
 \hspace*{8.5cm}
\resizebox*{13cm}{!}{\includegraphics{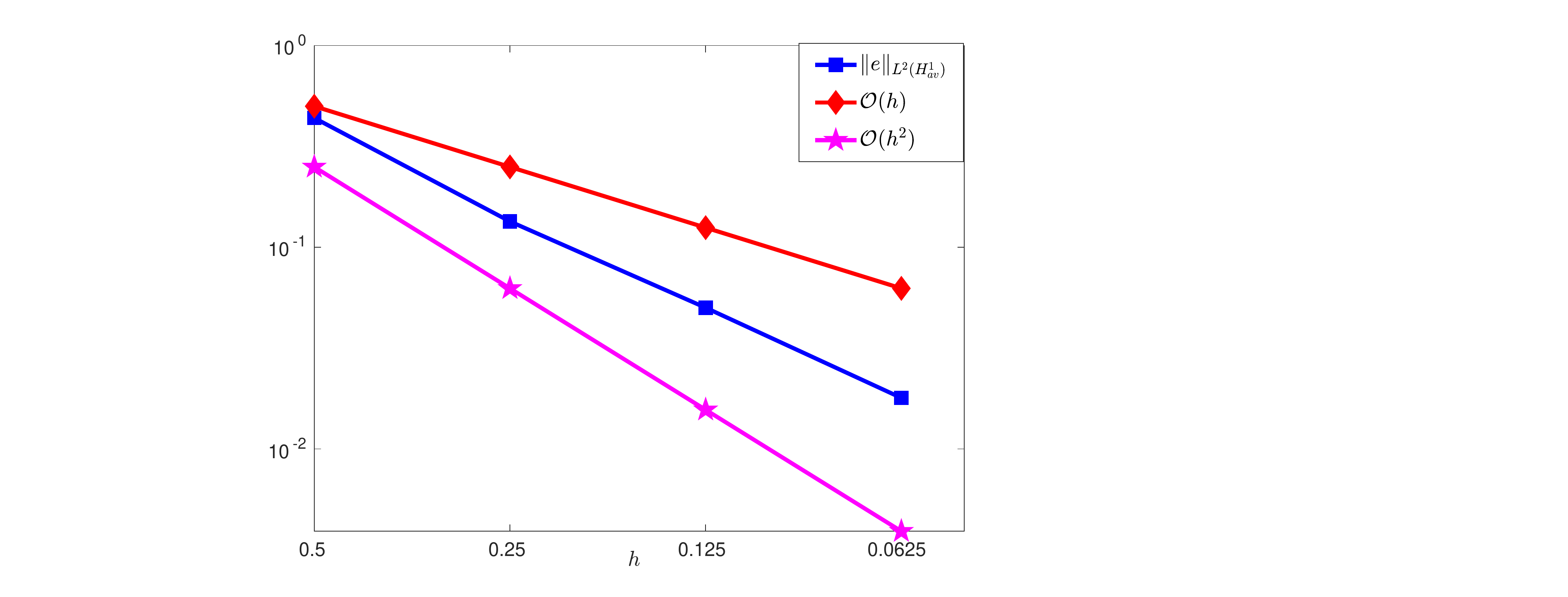}}%
\end{tabular}
} 
\caption{\label{fig:P23d} Errors for example 5.2 (3d)  for $P_2$ finite elements under simultaneous refinement in space and time ($\Delta t = h/10$) for $P_1$ finite elements.  \textbf{Top left}: $L^2$-norm at the end time. \textbf{Top right}: $L^2(L^2)$-norm. \textbf{Bottom}: $L^2(H_{av}^1)$-norm.}
\end{figure}

In Figure~\ref{fig:P23d}, we illustrate the errors under simultaneous refinement ($h_i=\frac{\Delta t_i}{10}$) for $P_2$ finite elements. We observe again (at least) second-order convergence in the $L^2$-norms, in agreement with Theorem~\ref{globalerrthm}. The convergence in the $L^2(H_{\text{av}}^1)$-norm lies between linear and quadratic convergence and decreases slightly for finer mesh sizes. This reduction in the convergence rate could be related to a violation of the CFL condition, which has been used in the error analysis. On the other hand, we did not find any stability issues in our computations.


\section{Concluding remarks}\label{Sec:6}

We have analysed a Crank-Nicolson variant of the implicitly extended Eulerian time-stepping scheme for the heat equation on time-dependent domains. Theoretically, stability and optimal-order convergence estimates were derived in the energy norm under the assumption of a parabolic CFL condition. 
In the numerical results, on the other hand, we did not observe any stability issue related to a violated CFL condition.
The three-dimensional results for second-order polynomials indicate that a violated CFL condition could result in a slightly reduced convergence order in the $L^2(H_{\text{av}}^1)$-norm.


To our knowledge this is the first work, in which an implicitly extended Eulerian time-stepping scheme is applied with a scheme that requires derivative information at different time steps. As mentioned in the introduction, this could be the basis for an analysis of a whole zoo of time-stepping schemes, such as the Fractional-step-$\theta$ method, implicit Runge-Kutta- or Adams-Bashforth schemes.
Moreover, we plan to apply the developed time-stepping scheme to flow problems on time-dependent domains and to fluid-structure interactions with large displacements, see e.g.~\cite{FreiDiss, BurmanFernandezFrei2020}.

%
%



\medskip
\bibliographystyle{plain}

\end{document}